\documentclass[11pt]{article}
\usepackage[T1]{fontenc}
\usepackage[utf8x]{inputenc}
\usepackage{color}
\usepackage{amstext}
\usepackage{amssymb}
\usepackage{esint}
\usepackage{authblk}
\usepackage[unicode=true,
 bookmarks=true,bookmarksnumbered=false,bookmarksopen=false,
 breaklinks=false,pdfborder={0 0 1},backref=false,colorlinks=false]
 {hyperref}

\makeatletter

\newcommand{\lyxmathsym}[1]{\ifmmode\begingroup\def\b@ld{bold}
  \text{\ifx\math@version\b@ld\bfseries\fi#1}\endgroup\else#1\fi}

\@ifundefined{date}{}{\date{}}

\def\Ba{{\cal B}}
\def\Ca{{\cal C}}
\def\Da{{\cal D}}
\def\Ea{{\cal E}}
\def\Fa{{\cal F}}
\def\Ga{{\cal G}}

\def\Ka{{\cal K}}
\def\La{{\cal L}}
\def\Ma{{\cal M}}

\def\Pa{{\cal P}}
\def\Qa{{\cal Q}}
\def\Ra{{\cal R}}
\def\Sa{{\cal S}}

\def\Ua{{\cal U}}
\def\Va{{\cal V}}
\def\Wa{{\cal W}}
\def\Xa{{\cal X}}

\def\QQ{\mathbb{Q}}
\def\PP{\mathbb{P}}
\def\EE{\mathbb{E}}
\def\NN{\mathbb{N}}

\def\RR{\mathbb{R}}

\newcommand{\point}{\mbox{\LARGE .}}

\def\blbx{\hbox{\vrule height 5pt width 5pt depth 0pt}\medskip}

\newcommand{\proof}{\textbf{Proof:}\hskip.3cm}

\newcommand{\cqfd}{\hfill\blbx \\}

\newtheorem{prop}{Proposition}[section]
\newtheorem{cor}[prop]{Corollary}
\newtheorem{defi}[prop]{Definition}

\newtheorem{lem}[prop]{Lemma}

\newtheorem{theo}[prop]{Theorem}
\newenvironment{rem}{\refstepcounter{prop}
{\bf Remark \theprop :\ }}{\par\me}

\newcommand{\bq}{\begin{eqnarray*}}
\newcommand{\eq}{\end{eqnarray*}}

\newcommand{\bqn}[1]{\begin{eqnarray}\label{#1}}
\newcommand{\eqn}{\end{eqnarray}}

\newcommand{\me}{\medskip}

\def\dessus#1#2{\mathord{\mathop{\kern 0pt #2}\limits^#1}}

\makeindex  

\makeatother

\begin{document}
\author[1]{Pierre Del Moral}
\author[2]{Pierre E. Jacob}
\author[3]{Anthony Lee}
\author[4]{Lawrence Murray}
\author[5]{Gareth W. Peters}

\affil[1]{INRIA Bordeaux Sud-Ouest and University of Bordeaux, France}
\affil[2]{National University of Singapore, Singapore}
\affil[3]{University of Warwick, UK}
\affil[4]{CSIRO Mathematics, Informatics and Statistics, Perth, Australia}
\affil[5]{University College London, UK}

\title{Feynman-Kac particle integration with geometric interacting jumps}


\maketitle
\begin{abstract}
This article is concerned with the design and analysis of discrete
time Feynman-Kac particle integration models with geometric interacting
jump processes. We analyze two general types of model, corresponding
to whether the reference process is in continuous or discrete time.
For the former, we consider discrete generation particle models defined
by arbitrarily fine time mesh approximations of the Feynman-Kac models
with continuous time path integrals. For the latter, we assume that the discrete process is observed at integer times and we design new approximation models with geometric interacting jumps in terms of
a sequence of intermediate time steps between the integers. In both
situations, we provide non asymptotic bias and variance theorems w.r.t.
the time step and the size of the system, yielding what appear to
be the first results of this type for this class of Feynman-Kac particle
integration models. We also discuss uniform convergence estimates
w.r.t. the time horizon. Our approach is based on an original semigroup
analysis with first order decompositions of the fluctuation errors.\\


\end{abstract}

\global\long\def\Ba{{\cal B}}
\global\long\def\Ca{{\cal C}}
\global\long\def\Da{{\cal D}}
\global\long\def\Ea{{\cal E}}
\global\long\def\Fa{{\cal F}}
\global\long\def\Ga{{\cal G}}
\global\long\def\Ka{{\cal K}}
\global\long\def\La{{\cal L}}
\global\long\def\Ma{{\cal M}}
\global\long\def\Pa{{\cal P}}
\global\long\def\Qa{{\cal Q}}
\global\long\def\Ra{{\cal R}}
\global\long\def\Sa{{\cal S}}
\global\long\def\Ua{{\cal U}}
\global\long\def\Va{{\cal V}}
\global\long\def\Wa{{\cal W}}
\global\long\def\Xa{{\cal X}}

\global\long\def\EE{\mathbb{E}}
\global\long\def\NN{\mathbb{N}}
\global\long\def\PP{\mathbb{P}}
\global\long\def\QQ{\mathbb{Q}}
\global\long\def\RR{\mathbb{R}}

\global\long\def\point{\mbox{\LARGE.}}

\section{Introduction}

Feynman-Kac formulae are central path integration mathematical models
in physics and probability theory. More precisely, these models and
their interacting particle interpretations have come to play a significant
role in applied probability, numerical physics, Bayesian statistics,
probabilistic machine learning, and engineering sciences. Applications
of these particle integration techniques are increasingly used to
solve a variety of complex problems in nonlinear filtering, data assimilation,
rare event sampling, hidden Markov chain parameter estimation, stochastic
control and financial mathematics. A detailed account of these functional
models and their application domains can be found in the series of
research books~\cite{cappe-2005,dm2004,dm2000a,doucet-2001,jefferies}
and, more recently, in~\cite{carmona,dhw2012,del2012introduction}.

In computational physics, these techniques are used for free energy
computations, specifically in estimating ground states of Schr\"odinger
operators. In this context, these particle models are often referred
as quantum or diffusion Monte Carlo methods~\cite{assaraf2000computing,assaraf2000diffusion,PhysRevB.85.115115,1996AmJPh..64..633K}.
We also refer the reader to the series of articles~\cite{dmsoft-2004,dm2003,rousset2007,rousset2010,rousset2006}.

In advanced signal processing, they are known as particle filters
or sequential Monte Carlo methods, and were introduced in three independent
works in the 90's~\cite{dm-96,gordon-93,kitagawa}. These stochastic
particle algorithms are now routinely used to compute sequentially
the flow of conditional distributions of the random states of a signal
process given some noisy and partial observations~\cite{cappe-2005,dm2004,ddj2006,dm2000a,doucet-2001,doucet-2011,kantas,kitagawa2}.
Feynman-Kac formulae and their particle interpretations are also commonly
used in financial mathematics to model option prices, futures prices
and sensitivity measures, and in insurance and risk models~\cite{carmona,fouque09,hugonnier,peters2011calibration,peters2007simulation,peters2009sequential}.
They are used in rare event analysis to model conditional distributions of stochastic processes evolving in a rare event regime~\cite{cerou-2006,fouque09,lezaud-2006}.

This article presents geometric interacting jump particle approximations
of Feynman-Kac path integrals. It also contains theoretical results
related to the practical implementation of these particle algorithms
for both discrete and continuous time integration problems. A key
result is the presentation of connections between the interacting
jump particle interpretations of the continuous time models and their
discrete time generation versions. This is motivated by the fact that
while the continuous time nature of these models is fundamental to
describing certain phenomena, the practical implementation of these
models on a computer requires a judicious choice of time discretization.
Conversely, as shown in section 2.1 in~\cite{dm2007}, a discrete
time Feynman-Kac model can be encapsulated within a continuous time
framework by considering stochastic processes only varying on integer
times. Continuous time Feynman-Kac particle models are based on exponential
interacting jumps~\cite{dh-2010,dm2003,dm2000a,eberle1,eberle2,eberle3,rousset2006}, while
their discrete time versions are based on geometric type jumps~\cite{dm2004,dhw2012,dkm-2001}.
From a computational perspective, the exponential type interacting
jumps thus need to be approximated by geometric type jumps. Incidentally,
some of these geometric type interacting jump particle algorithms
are better suited to implementation in a parallel computing environment
(see section \ref{remark:parallelcomputation}).

Surprisingly, little attention has been paid to analyze the connections
between exponential and geometric type jump particle models. There
are references dealing with these two models separately~\cite{crisan-99,dmj1,dmj2,dm2000b,dm2007,dm2003,rousset2006},
but none provide a convergence analysis between the two. In this paper
we initiate this study with a non asymptotic bias and variance analysis
w.r.t. the time step parameter and the size of the particle population
scheme. Special attention is paid to the stochastic modeling of these
interacting jump processes, and to a stochastic perturbation analysis
of these particle models w.r.t. local sampling random fields\textcolor{red}{}.

We conclude this section with basic notation used in the article.
We let $\Ba_{b}(E)$ be the Banach space of all bounded Borel functions
$f$ on some Polish%
\footnote{i.e. homeomorphic to a complete separable metric space%
} state space $E$ equipped with a Borel $\sigma$-field $\Ea$, equipped
with the uniform norm $\|f\|=\sup_{x\in E}|f(x)|$. We denote by $\mbox{{\rm osc}}(f):=\sup_{x,y}{|f(x)-f(y)|}$
the oscillation of a function $f\in\Ba_{b}(E)$. We let $\mu(f)=\int f(x)\mu(dx)$
be the Lebesgue integral of a function $f\in\Ba_{b}(E)$ with respect
to a finite signed measure $\mu$ on $E$. We also equip the set $\Ma(E)$
of finite signed measures $\mu$ with the total variation norm $\|\mu\|_{\rm tv}=\sup_{}{|\mu(f)|}$,
where the supremum is taken over all functions $f\in\Ba_{b}(E)$ with
$\mbox{{\rm osc}}(f)\leq1$. We let $\Pa(E)\subset\Ma(E)$ be the
subset of all probability measures. We recall that any bounded integral
operator $Q$ on $E$ is an operator $Q$ from $\Ba_{b}(E)$ into
itself defined by $Q(f)(x)=\int Q(x,dy)f(y)$, for some measure $Q(x,\point)$,
indexed by $x\in E$, and we set $\|Q\|_{\rm tv}=\sup_{x\in E}\|Q(x,\point)\|_{\rm tv}$.
These operators generate a dual operator $\mu\mapsto\mu Q$ on the
set of finite signed measures defined by $(\mu Q)(f)=\mu(Q(f))$.
A Markov kernel is a positive and bounded integral operator $Q$ s.t.
$Q(1)=1$. The Dobrushin contraction coefficient of a Markov kernel
$Q$ is defined by $\beta(Q):=\sup{\mbox{{\rm osc}}(Q(f))}$, where
the supremum is taken over all functions $f\in\Ba_{b}(E)$ s.t. $\mbox{{\rm osc}}(f)\leq1$.
Given some positive potential function $G$ on $E$, we denote by
$\Psi_{G}$ the Boltzmann-Gibbs transformation $\mu\in\Pa(E)\mapsto\Psi_{G}(\mu)\in\Pa(E)$
defined by $\Psi_{G}(\mu)(f)=\mu(fG)/\mu(G)$.

\section{Description of the models}

\subsection{Feynman-Kac models}

We consider an $E$-valued Markov process $\Xa_{t}$, $t\in\RR_{+}=[0,\infty[$
defined on a standard filtered probability space $(\Omega,\Fa=(\Fa_{t})_{t\in\RR_{+}},\PP)$.
The set $\Omega=D(\RR_{+},E)$ represents the space of càdlàg paths
equipped with the Skorokhod topology which turn it into a Polish space.
A point $\omega\in\Omega$ represents a sample path of the canonical
process $\Xa_{t}(\omega)=\omega_{t}$. We also let $\Fa_{t}^{X}=\sigma(\Xa_{s},~s\leq t)$
and $\PP$ be the sigma-field and probability measure of the process
$(\Xa_{t})_{t\in\RR_{+}}$. Finally, we also consider the $\PP$-augmentation
$\Fa_{t}$ of $\Fa_{t}^{X}$ so that the resulting filtration satisfies
the usual conditions of right continuity and completion by $\PP$-negligible
sets (see for instance~\cite{karat,protter}, and the references
therein). We also consider a time inhomogeneous bounded Borel function
$\Va_{t}$ on $E$.

We let $\QQ_{t}$ and $\Lambda_{t}$ be the Feynman-Kac measures on
$\Omega_{t}:=D([0,t],E)$ defined for any bounded measurable function
$f$ on $\Omega_{t}$, by the following formulae 
\begin{equation}
\QQ_{t}(f):=\Lambda_{t}(f)/\Lambda_{t}(1)\ \quad\mbox{{\rm with}}\quad\Lambda_{t}(f)=\EE\left(f(\Xa_{[0,t]})~\exp{\left(\int_{0}^{t}\Va_{s}(\Xa_{s})ds\right)}\right)\label{FKC}
\end{equation}
and we let $\nu_{t}$ and $\mu_{t}$, respectively, be the $t$-marginals
of $\Lambda_{t}$ and $\QQ_{t}$.

We consider the mesh sequence $t_{k}=k/m$, $k\geq0$, with time step
$h=t_{n}-t_{n-1}=1/m$ associated with some integer $m\geq1$, and
we let $\QQ_{t_{n}}^{(m)}$ and $\Lambda_{t_{n}}^{(m)}$ be the Feynman-Kac
measures on $\Omega_{t_{n}}$ defined for any bounded measurable function
$f$ on $\Omega_{t_{n}}$, by the following formulae 
\begin{equation}
\QQ_{t_{n}}^{(m)}(f):=\Lambda_{t_{n}}^{(m)}(f)/\Lambda_{t_{n}}^{(m)}(1)\quad\mbox{{\rm with}}\quad\Lambda_{t_{n}}^{(m)}(f)=\EE\left(f(\Xa_{[0,t_{n}]})~\prod_{0\leq p<n}e^{\Va_{t_{p}}(\Xa_{t_{p}})/m}~\right).\label{FKCD}
\end{equation}
We also denote by $\nu_{t_{n}}^{(m)}$ and $\mu_{t_{n}}^{(m)}$, respectively,
the $t_{n}$-marginal of $\Lambda_{t_{n}}^{(m)}$ and $\QQ_{t_{n}}^{(m)}$.
\begin{itemize}
\item \textbf{Case (D) :} {\em We have $\Xa_{t}=X_{\lfloor t\rfloor}$
and $\Va_{t}=\log{G_{\lfloor t\rfloor}}$, where $X_{n}$, $n\in\NN$
is an $E$-valued Markov chain, and $G_{n}$ are Borel positive functions
s.t. $\log{G_{n}}$ is bounded. }

In this case, the marginal $\nu_{n}=\gamma_{n}$ and $\mu_{n}=\eta_{n}$
of the Feynman-Kac measures of $\Lambda_{t}$ and $\QQ_{t}$ on integer
times $t=n$ are given for any $f\in\Ba_{b}(E)$ by the formula 
\begin{equation}
\eta_{n}(f)=\gamma_{n}(f)/\gamma_{n}(1)\quad\mbox{{\rm with}}\quad\gamma_{n}(f):=\EE\left(f(X_{n})~\prod_{0\leq p<n}G_{p}(X_{p})\right).\label{FKD}
\end{equation}

\item \textbf{Case (C) :} \label{casC-ref} {\em The process $\Xa_{t}$
is a continuous time Markov process with infinitesimal generators
$L_{t}~:~D(L)\rightarrow D(L)$ defined on some common domain of functions
$D(L)$, and $\Va\in\Ca^{1}(\RR_{+},D(L))$. The set $D(L)$ is a
sub-algebra of the Banach space $\Ba_{b}(E)$ generating the Borel
$\sigma$-field $\Ea$, and for any measurable function $\Ua:t\in\RR_{+}\mapsto\Ua_{t}\in\Da(L)$
the Feynman-Kac semigroup $\Qa_{s,t}$, $s\leq t$, defined by 
\[
\Qa_{s,t}(f)(x)=\EE\left(f(\Xa_{t})~\exp{\left(\int_{s}^{t}\Ua_{r}(\Xa_{r})dr\right)}~\left|~\Xa_{s}=x\right.\right)
\]
leaves $D(L)$ invariant; that is we have that $\Qa_{s,t}(D(L))\subset D(L)$.
For any $s\leq t$, the mappings $r\in[0,s]\mapsto L_{r}(Q_{s,t}(f)^{2})$
and $r\in[0,s]\mapsto L_{r}^{2}(Q_{s,t}(f)^{2})\in\Ca^{1}([0,s],D(L))$,
and their norm as well the norm of the first order derivatives only
depend on $(\Xa_{s})_{s\leq t}$ and on the norms of the functions
$(\Ua_{s})_{s\leq t}$ and their derivatives. }\label{CaseC} 
\end{itemize}
The regularity conditions stated in \textbf{(C)} correspond to time
inhomogeneous versions of those introduced in~\cite{dm2003}. They
hold for pure jump processes with bounded jump rates with $D(L)=\Ba_{b}(E)$,
or for Euclidean diffusions on $E=\RR^{d}$ with regular and Lipschitz
coefficients by taking $D(L)$ as the set of $\Ca^{\infty}$-functions
with derivatives decreasing at infinity faster that any polynomial
function. These regularity conditions allow the use of most of the
principal theorems of stochastic differential calculus, e.g. the ``carré
du champ'', or square field, operator that characterizes the predictable
quadratic variations of the martingales that appear in Ito's formulae.
These regularity conditions can probably be relaxed using the extended
setup developed in~\cite{dm2000a}.\label{page-ex}

We have already mentioned that the particle interpretations associated
with the continuous time models (\ref{FKC}) are defined in terms
of interacting jump particle systems~\cite{dm2000a,dm2000b,dm2003,dm2007}.
The implementation of these continuous time particle algorithms is
clearly impractical and we therefore resort to the geometric interacting
processes associated with the $m$-approximation models defined in
(\ref{FKCD}). These discrete generation interacting jumps models
provide new and different types of adaptive resampling procedures,
which differ from those discussed in the articles~\cite{ddj2006,ddj-2008},
and the references therein.

\subsection{Mean field particle models}

In this section, we provide a brief description of the geometric type
interacting jump particle models associated with the $m$-approximation
Feynman-Kac model defined in (\ref{FKCD}). First, if we define 
\[
\Ma_{t_{n},t_{n+1}}(x,dy)=\PP\left(\Xa_{t_{n+1}}\in dy~|~\Xa_{t_{n}}=x\right)\quad\mbox{{\rm and}}\quad\Ga_{t_{n}}=\exp{\left(\Va_{t_{n}}/m\right)},
\]
then it is well known that $\mu_{t_{n}}^{(m)}$ satisfies the following
evolution equation 
\begin{equation}
\mu_{t_{n+1}}^{(m)}=\Psi_{\Ga_{t_{n}}}(\mu_{t_{n}}^{(m)})\Ma_{t_{n},t_{n+1}}.\label{GENC}
\end{equation}
Further details on the derivation of these evolution equations can
be found in~\cite{dm2004,dm2000a,dhw2012}. The particle interpretation
of this model depends on the interpretation of the Boltzmann-Gibbs
transformation in terms of a Markov transport equation 
\begin{equation}
\Psi_{\Ga_{t_{n}}}(\mu)=\mu\Sa_{t_{n},\mu}\label{MT}
\end{equation}
for some Markov transitions $\Sa_{t_{n},\mu}$, that depend on the
time parameter $t_{n}$ and on the measure $\mu$. The choice of these
Markov operators is not unique; we refer to~\cite{dm2004} for a
more thorough discussion of these models. In this article, we consider
an abstract general model, and illustrate our study with the following
three classes of models. 
\begin{itemize}
\item \textbf{Case 1 :}\label{case1-ref} We have $\Va_{t}=-\Ua_{t}$, for
some non negative and bounded function $\Ua_{t}$. In this situation,
(\ref{MT}) is satisfied by the Markov transition 
\[
\Sa_{t_{n},\mu}(x,dy):=e^{-\Ua_{t_{n}}(x)/m}~\delta_{x}(dy)+\left(1-e^{-\Ua_{t_{n}}(x)/m}\right)~\Psi_{e^{-\Ua_{t_{n}}/m}}(\mu)(dy).
\]

\item \textbf{Case 2 :} The function $\Va_{t}$ is non negative. In this
situation, (\ref{MT}) is satisfied by the Markov transition 
\[
\Sa_{t_{n},\mu}(x,dy):=\frac{1}{\mu\left(e^{\Va_{t_{n}}/m}\right)}~\delta_{x}(dy)+\left(1-\frac{1}{\mu\left(e^{\Va_{t_{n}}/m}\right)}\right)~\Psi_{\left(e^{\Va_{t_{n}}/m}-1\right)}(\mu)(dy).
\]

\item \textbf{Case 3:} The Markov transport equation (\ref{MT}) is satisfied
by the Markov transition 
\[
\Sa_{t_{n},\mu}(x,dy):=\left(1-a_{t_{n},\mu}(x)\right)~\delta_{x}(dy)+a_{t_{n},\mu}(x)~\Psi_{\left(e^{\Va_{t_{n}}/m}-e^{\Va_{t_{n}}(x)/m}\right)_{+}}(\mu)(dy)
\]
with the rejection rate $a_{t_{n},\mu}(x):={\mu\left(\left(e^{\Va_{t_{n}}/m}-e^{\Va_{t_{n}}(x)/m}\right)_{+}\right)}/{\mu\left(e^{\Va_{t_{n}}/m}\right)}\in[0,1]$ 
\end{itemize}
In these three cases we have the following first order expansion 
\begin{equation}
\Sa_{t_{n},\mu}=Id+\frac{1}{m}~\widehat{L}_{t_{n},\mu}+~\frac{1}{m^{2}}~\widehat{R}_{t_{n},\mu}\label{1st-exp-ref}
\end{equation}
with some jump type generator $\widehat{L}_{t_{n},\mu}$ and some
integral operator $\widehat{R}_{t_{n},\mu_{t_{n}}^{(m)}}$ s.t. $\sup_{}{\left\Vert \widehat{R}_{t_{n},\mu}\right\Vert _{{\rm tv}}}<\infty$,
where the supremum is taken over all $m\geq1$ and $\mu\in\Pa(E)$.
The jump generators $\widehat{L}_{t_{n},\mu}$ corresponding to the
three cases presented above are described respectively in (\ref{gen1-ref}),
(\ref{gen2-ref}), and (\ref{gen3-ref}). The proofs of these expansions
is rather elementary, and they are housed in the appendix, on page~\pageref{decomp-a}.

In addition, whenever (\ref{MT}) is satisfied, we have the evolution
equation 
\begin{equation}
\mu_{t_{n+1}}^{(m)}=\mu_{t_{n}}^{(m)}\Ka_{n+1,\mu_{t_{n}}^{(m)}}\quad\mbox{{\rm with the Markov kernels}}\quad\Ka_{t_{n},t_{n+1},\mu}=\Sa_{t_{n},\mu}\Ma_{t_{n},t_{n+1}}.\label{WEAKD}
\end{equation}
The mean field $N$-particle model $\xi_{t_{n}}:=\left(\xi_{t_{n}}^{i}\right)_{1\leq i\leq N}$
associated with the evolution equation (\ref{WEAKD}) is a Markov
process in $E^{N}$ with elementary transitions given by 
\begin{equation}
\PP\left(\xi_{t_{n+1}}\in dx~|~ \xi_{t_{n}}\right)=\prod_{1\leq i\leq N}\Ka_{t_{n},t_{n+1},\mu_{t_{n}}^{N}}(\xi_{t_{n}}^{i},dx^{i})\quad\mbox{{\rm with}}\quad\mu_{t_{n}}^{N}=\frac{1}{N}\sum_{1\leq i\leq N}\delta_{\xi_{t_{n}}^{i}},\label{MF-ref}
\end{equation}
where $dx=dx^{1}\ldots dx^{N}$ stands for an infinitesimal neighborhood
of the point $x=(x^{i})_{1\leq i\leq N}\in E^{N}$.

\section{Statement of the main results}

Our first main result relates the Feynman-Kac models (\ref{FKC})
and their $m$-approximation measures (\ref{FKCD}) in case \textbf{(D)}
and \textbf{(C)}.

\begin{theo}\label{theo1}

In case \textbf{(D)}, we have 
\[
\nu_{n}^{(m)}=\nu_{n}=\gamma_{n}\quad\mbox{and}\quad\mu_{n}^{(m)}=\mu_{n}=\eta_{n}
\]
with the Feynman-Kac measures $\gamma_{n}$ and $\eta_{n}$ defined
in (\ref{FKD}).

In case \textbf{(C)}, we have the first order decomposition 
\[
\Lambda_{t_{n}}^{(m)}=\Lambda_{t_{n}}+\frac{1}{m}~r_{m,t_{n}}\quad\mbox{and}\quad\QQ_{t_{n}}^{(m)}=\QQ_{t_{n}}+\frac{1}{m}~\overline{r}_{m,t_{n}}
\]
with some remainder signed measures $r_{m,t_{n}},\overline{r}_{m,t_{n}}$
s.t. $\sup_{m\geq1}{\left[\left\Vert \overline{r}_{m,t_{n}}\right\Vert _{\rm tv}\vee\left\Vert r_{m,t_{n}}\right\Vert _{\rm tv}\right]}<\infty$.

\end{theo}

The proof of the theorem is rather technical and it is postponed to
the appendix.

The first assertion of theorem~\ref{theo1} allows us to turn a discrete
time Feynman-Kac model (\ref{FKD}) into a continuous time model (\ref{FKC}).
To be more precise, we have that $\nu_{t_{n}}^{(m)}=\nu_{t_{p}}^{(m)}Q_{t_{p},t_{n}}^{(m)}$
with the Feynman-Kac semigroup 
\[
Q_{t_{p},t_{n}}^{(m)}(f)(x):=\EE\left(f(\Xa_{t_{n}})~\prod_{p\leq q<n}e^{\Va_{t_{q}}(\Xa_{t_{q}})/m}~\left|~\Xa_{t_{p}}=x\right.\right)
\]
in case \textbf{(D)}, for integer times $(t_{p},t_{n})=(km,nm)$,
with $k\leq n$, we also have that $\gamma_{n}=\gamma_{k}Q_{k,n}$
with the Feynman-Kac semigroup 
\[
Q_{k,n}(f)(x):=\EE\left(f(X_{n})~\prod_{k\leq l<n}G_{l}(X_{l})~\left|~X_{k}=x\right.\right)=Q_{k,n}^{(m)}(f)(x).
\]
Thus, the normalized Markov kernels $P_{k,n}^{(m)}(f):=Q_{k,n}^{(m)}(f)/Q_{k,n}^{(m)}(1)$
also coincide with the Markov kernels $P_{k,n}(f):=Q_{k,n}(f)/Q_{k,n}(1)$.
In addition, for any $k\geq0$ and $r<m$, we also have the semigroup
formulae 
\begin{equation}
Q_{k,n}^{(m)}(f)(x)=G_{k}(x)^{r/m}~Q_{{k+r/m},n}^{(m)}(f)(x)\quad\mbox{{\rm and}}\quad P_{k,n}^{(m)}=P_{k+r/m,n}^{(m)}=P_{k,n}.\label{sg-formula}
\end{equation}

We prove the l.h.s. assertion using the fact that for any $n\geq0$
and any $p=km+r$, with $k\geq0$ and $r<m$, we have $t_{p}=k+r/m$
and 
\begin{eqnarray*}
Q_{k,n}^{(m)}(f)(x) & = & G_{k}(x)^{r/m}\times\EE\left(f(\Xa_{t_{nm}})~\prod_{k+r/m\leq t_{q}<n}e^{V_{t_{q}}(\Xa_{t_{q}})/m}~\left|~\Xa_{{k+r/m}}=x\right.\right).
\end{eqnarray*}

For a Feynman-Kac measure (\ref{FKC}) associated with a continuous
diffusion style process $\Xa_{t}$, it is important to observe that
the l.h.s. measure in the $m$-approximation model (\ref{FKCD}),
as defined on a time mesh sequence, can be thought of as a time discretization
of the exponential path integrals in the continuous time model (\ref{FKC}).
Nevertheless, the elementary Markov transitions of the Markov chain
$(\Xa_{t_{n}})_{n\geq0}$ are generally unknown. To get some feasible
Monte Carlo approximation scheme, we need a dedicated technique to
sample the transitions of this chain. One natural strategy is to replace
in (\ref{FKCD}), the reference Markov chain $(\Xa_{t_{n}})_{n\geq0}$
by the Markov chain $(\hat{\Xa}_{t_{n}})_{n\geq0}$ associated with
some Euler type discretization model with time step $\Delta t=1/m$.
The stochastic analysis of these models is discussed in some detail
in the articles~\cite{dmj1,dmj2,dh-2010}, including first order
expansions in terms of the size of the time mesh sequence.

Our second main result is the following non asymptotic bias and variance
theorem for the $N$-approximation mean field model introduced in
(\ref{MF-ref}).

\begin{theo}\label{theo-bv1-intro} We assume that the Markov transport
equation (\ref{MT}) is satisfied for Markov transitions $\Sa_{t_{n},\mu}$
also satisfying the first order decomposition (\ref{1st-exp-ref}).

In case \textbf{(C)}, for any function $f\in D(L)$, and any $N\geq m\geq1$
we have the non asymptotic bias and variance estimates 
\[
\left|\EE\left(\mu_{t_{n}}^{N}(f)\right)-\mu_{t_{n}}(f)\right|\leq c_{t_{n}}(f)~\left[\frac{1}{N}+\frac{1}{m}\right]
\]
and 
\[
\EE\left(\left[\mu_{t_{n}}^{N}(f)-\mu_{t_{n}}(f)\right]^{2}\right)\leq c_{t_{n}}(f)~\left[\frac{1}{N}+\frac{1}{m^{2}}\right]
\]
for some finite constant $c_{t_{n}}(f)<\infty$ that only depends
on $t_{n}$ and on $f$.

In case \textbf{(D)}, for any $f\in\Ba_{b}(E)$ s.t. $\mbox{{\rm osc}}(f)\leq1$,
and for any $N\geq m\geq1$ we have the non asymptotic bias and variance estimates
\[
N~\left|\EE\left(\mu_{n}^{N}(f)\right)-\eta_{n}(f)\right|\leq a(n)\quad\mbox{and}\quad N~\EE\left(\left[\mu_{t_{n}}^{N}(f)-\mu_{t_{n}}(f)\right]^{2}\right)\leq a(n)~\left(1+\frac{1}{N}~a(n)\right)
\]
for some some constant 
\[
a(n)\leq c~{\displaystyle \sum_{0\leq k<n}g_{k,n}^{3}g_{k,k+1}^{3}\left(\|\log{G_{k}}\|\vee1\right)^{2}\beta\left(P_{k,n}\right)\quad\mbox{with}\quad g_{k,n}:=\sup_{x,y}{Q_{k,n}(1)(x)/Q_{k,n}(1)(y})}.
\]

\end{theo}

Under appropriate regularity conditions on the Feynman-Kac model,
we can prove that the constant $a(n)$ is uniformly bounded w.r.t.
the time parameter; that is we have that $\sup_{n\geq0}a(n)<\infty$.
For a detailed discussion of these uniform convergence properties
w.r.t. the time parameter, we refer the reader to the book~\cite{dm2004},
and the more recent article~\cite{dhw2012}. To be more precise,
we let $\Phi_{k,l}(\eta_{k})=\eta_{l}$ be the Feynman-Kac semigroup
associated with the flow of measures $\eta_{k}$. In this notation,
by proposition 2.3 in~\cite{dm-Toulouse} we have that the Dobrushin contraction coefficient of the Markov kernel $Q_{k,n}(f)/Q_{k,n}(1)$ is given by
\[
\beta(P_{k,n})=\sup_{\mu_{1},\mu_{2}\in\Pa(E)}\left\Vert \Phi_{k,n}(\mu_{1})-\Phi_{k,n}(\mu_{2})\right\Vert _{\rm tv}.
\]
On the other hand, we also have that 
\[
Q_{k,n}(1)(x)=\prod_{k\leq l<n}\Phi_{k,l}(\delta_{x})(G_{l})\Rightarrow\log{\frac{Q_{k,n}(1)(x)}{Q_{k,n}(1)(y)}}=\sum_{k\leq l<n}\left(\log{\Phi_{k,l}(\delta_{x})(G_{l})}-\log{\Phi_{k,l}(\delta_{y})(G_{l})}\right).
\]
Using the fact that $\log{a}-\log{b}=\int_{0}^{1}\frac{(a-b)}{ta+(1-t)b}dt$,
we find that 
\[
\log{\frac{Q_{k,n}(1)(x)}{Q_{k,n}(1)(y)}}=\sum_{k\leq l<n}\int_{0}^{1}\frac{\left[\Phi_{k,l}(\delta_{x})(G_{l})-\Phi_{k,l}(\delta_{x})(G_{l})\right]}{t\Phi_{k,l}(\delta_{x})(G_{l})+(1-t)\Phi_{k,l}(\delta_{y})(G_{l})}~dt.
\]
Assuming that for any $l$ and $x$, and 
\begin{equation}
c_{1}\leq G_{l}(x)\leq c_{2}\quad\mbox{{\rm and}}\quad\sup_{\mu_{1},\mu_{2}\in\Pa(E)}\left\Vert \Phi_{k,l}(\mu_{1})-\Phi_{k,l}(\mu_{2})\right\Vert _{\rm tv}\leq c_{3}~e^{-c_{4}(k-l)}\label{stab-condition}
\end{equation}
for some positive and bounded constants $c_{i}$, $1\leq i\leq4$,
we find that 
\[
\beta(P_{k,n})\leq c_{3}~e^{-c_{4}(k-n)}\quad\mbox{{\rm and}}\quad\log{g_{k,n}}\leq2(c_{2}c_{3}/c_{1})\left(\sum_{k\leq l<n}e^{-c_{4}(k-l)}\right)\leq2(c_{2}c_{3})/(c_{1}(1-e^{-c_{4}})).
\]
This clearly implies that $(\ref{stab-condition})\Rightarrow\sup_{n\geq0}a(n)<\infty$.

For instance, it was proven in \cite{dm2000a,dmg2001} that condition
(\ref{stab-condition}) is met for time homogeneous models as soon
as the Markov transition $M$ of the Markov chain $X_{n}$ satisfies
the following mixing condition 
\[
\exists m\geq1,~\exists\rho>0~:~\forall x,y\in E\quad M^{m}(x,\point)\geq\rho~M^{m}(y,\point).
\]
It is well known that this condition is satisfied for any aperiodic
and irreducible Markov chains on finite state spaces, as well as for
bi-Laplace exponential transitions associated with a bounded drift
function, and for Gaussian transitions with a mean drift function
that is constant outside some compact domain.

The remainder of the article is organized as follows:

Section~\ref{sec-cont-mod} is concerned with continuous time particle
interpretations of the Feynman-Kac models (\ref{FKC}). By the representation
theorem~\ref{theo1}, these schemes also provide a continuous time
particle interpretation of the discrete time models (\ref{FKD}) without
further work. In section~\ref{mckean-sec}, we present the McKean
interpretation of the Feynman-Kac models in terms of a time inhomogeneous
Markov process whose generator depends on the distribution of the
random states. The choice of these McKean models is not unique. We
discuss the three interpretation models corresponding to the three
selection type transitions presented on page~\pageref{case1-ref}.
The mean field particle interpretation of these McKean models are
discussed in section~\ref{mean-fieldC}.

Of course, even for discrete time models (\ref{FKD}) these continuous
time particle interpretations are based on continuous time interacting
jump models and they cannot be used in practice without an additional
level of approximation. In this context, when using an Euler type
approximation these exponential interacting jumps are replaced by
geometric type recycling clocks. These interacting geometric jumps
particle models are discussed in section~\ref{sec-disc-mod}, which
is dedicated to the discrete time particle interpretations of the
Feynman-Kac models presented in (\ref{FKCD}). In section~\ref{mckean-secD},
we discuss the McKean interpretation of the Feynman-Kac models in
terms of a time inhomogeneous Markov chain model whose elementary
transitions depends on the distribution of the random states. Again,
the choice of these McKean models is not unique. We discuss the three
interpretation models corresponding to the three cases presented on
page~\pageref{case1-ref}. The mean field particle interpretation
of these McKean models are discussed on page~\pageref{TWO}.

Once again, using the representation formulae (\ref{FKCD}) we emphasize
that these schemes also provide a discrete generation particle interpretation
of the discrete time models (\ref{FKD}). In contrast to standard
discrete generation particle models associated with (\ref{FKD}),
these particle schemes are defined on a refined time mesh sequence
between integers. This time mesh sequence can be interpreted as a
time dilation. Between two integers, the particle evolution undergoes
an additional series of intermediate time evolution steps. In each
of these time steps, a dedicated Bernoulli acceptance-rejection trial
coupled with a recycling scheme is performed. As the time step decreases
to $0$, the resulting geometric interacting jump processes converge
to the exponential interacting jump processes associated with the
continuous time particle model. The final section, section~\ref{first-order-sect},
is mainly concerned with the proof of theorem~\ref{theo-bv1-intro}.

\section{Continuous time models}

\label{sec-cont-mod}

\subsection{Feynman-Kac semigroups}

In case \textbf{(C)} the semigroup of the flow of non negative measures
$\nu_{t}$ is given for any $s\leq t$ by the following formulae $\nu_{t}=\nu_{s}Q_{s,t}$,
with the Feynman-Kac semigroup $Q_{s,t}$ defined for any $f\in\Ba(E)$
by 
\[
Q_{s,t}(f)(x)=\EE\left(f(\Xa_{t})~\exp{\left\{ \int_{s}^{t}\Va_{s}(\Xa_{s})~ds\right\} }~|~ X_{s}=x\right).
\]
This yields $\mu_{t}=\Phi_{s,t}(\mu_{s})$, with the nonlinear transformation
$\Phi_{s,t}$ on the set of probability measures defined for any $f\in\Ba(E)$
by 
\[
\Phi_{s,t}(\mu_{s})(f):=\mu_{s}(Q_{s,t}(f))/\mu_{s}(Q_{s,t}(1)).
\]
Using some stochastic calculus manipulations, we readily prove that
$\mu_{t}$ satisfies the following integro-differential equation 
\begin{equation}
\frac{d}{dt}\mu_{t}(f)=\mu_{t}(L_{t}(f))+\mu_{t}(\Va_{t}f)-\mu_{t}(\Va_{t})\mu_{t}(f)\label{GEN}
\end{equation}
for any function $f\in D(L)$. Further details on the derivation of
these evolution equations can be found in the articles~\cite{dm2003,dm2000b}.
The particle interpretation of this model depends on the interpretation
of the correlation term in the r.h.s. of (\ref{GEN}) in terms of
a jump type generator. The choice of these generators is not unique.
Next, we discuss three important classes of models. These three situations
are the continuous time versions of the three cases discussed on page~\pageref{case1-ref}. 
\begin{itemize}
\item \textbf{Case 1 :} We assume that $\Va_{t}=-\Ua_{t}$, for some non
negative function $\Ua_{t}$. In this situation, we have the formula
\[
\mu_{t}(\Va_{t}f)-\mu_{t}(\Va_{t})\mu_{t}(f)=\mu_{t}(\Ua_{t}\left[\mu_{t}(f)-f\right])=\mu_{t}\left(\widehat{L}_{t,\mu_{t}}(f)\right)
\]
with the interacting jump generator 
\begin{equation}
\widehat{L}_{t,\mu_{t}}(f)(x)=\Ua_{t}(x)~\int~[f(y)-f(x)]~\mu_{t}(dy).\label{gen1-ref}
\end{equation}

\item \textbf{Case 2 :} When $\Va_{t}$ is a positive function, then we
have the formula 
\[
\mu_{t}(\Va_{t}f)-\mu_{t}(\Va_{t})\mu_{t}(f)=\mu_{t}\left(\widehat{L}_{t,\mu_{t}}(f)\right)
\]
with the interacting jump generator 
\begin{equation}
\widehat{L}_{t,\mu_{t}}(f)(x)=\int~[f(y)-f(x)]~\Va_{t}(y)~\mu_{t}(dy)=\mu_{t}(\Va_{t})~\int~[f(y)-f(x)]~\Psi_{\Va_{t}}(\mu_{t})(dy).\label{gen2-ref}
\end{equation}

\item \textbf{Case 3 :} For any bounded potential function $\Va_{t}$ we
have 
\[
\mu_{t}(\Va_{t}f)-\mu_{t}(\Va_{t})\mu_{t}(f)=\int(f(y)-f(x))~(\Va_{t}(y)-\Va_{t}(x))_{+}~\mu_{t}(dx)~\mu_{t}(dy)=\mu_{t}\left(\widehat{L}_{t,\mu_{t}}(f)\right)
\]
with $a_{+}=a\vee0$, and with the interacting jump generator 
\begin{equation}
\widehat{L}_{t,\mu_{t}}(f)(x)=\int~[f(y)-f(x)]~(\Va_{t}(y)-\Va_{t}(x))_{+}~\mu_{t}(dy).\label{gen3-ref}
\end{equation}

\end{itemize}

\subsection{McKean interpretation models}

\label{mckean-sec}

In the three cases discussed above, for any test functions $f\in D(L)$
we have the evolution equation 
\begin{equation}
\frac{d}{dt}\mu_{t}(f)=\mu_{t}(L_{t,\mu_{t}}(f))\quad\mbox{{\rm with}}\quad L_{t,\mu_{t}}:=L_{t}+\widehat{L}_{t,\mu_{t}}.\label{WEAK}
\end{equation}
These integro-differential equations can be interpreted as the evolution
of the laws, given by $\mbox{{\rm Law}}(\overline{X}_{t})=\mu_{t}$,
of a time inhomogeneous Markov process $\overline{X}_{t}$ with infinitesimal
generators $L_{t,\mu_{t}}$ that depend on the distribution of the
random states at the previous time increment. This probabilistic model is called the McKean interpretation
of the evolution equation (\ref{WEAK}) in terms of a time inhomogeneous
Markov process. In this framework, using Ito's formula for any test
function $f\in\Ca^{1}([0,\infty[,\Da(L))$, we have 
\begin{equation}
df_{t}(\overline{\Xa}_{t})=\left(\frac{\partial}{\partial t}+L_{t,\mu_{t}}\right)(f_{t})(\overline{\Xa}_{t})+d\overline{M}_{t}(f)\label{ItoC}
\end{equation}
with a martingale term $\overline{M}_{t}(f)$ with predictable angle
bracket 
\[
d\langle\overline{M}(f)\rangle_{t}=\Gamma_{L_{t,\mu_{t}}}(f_{t},f_{t})(\overline{\Xa}_{t})dt.
\]
Using the r.h.s. description of $L_{t,\mu_{t}}$ in (\ref{WEAK}),
for any $f\in D(L)$ we notice that 
\[
\Gamma_{L_{t,\mu_{t}}}(f,f)=\Gamma_{L_{t}}(f,f)+\Gamma_{\widehat{L}_{t,\mu_{t}}}(f,f).
\]

Next, we provide a description of this Markov process in the three
cases discussed above. 
\begin{itemize}
\item \textbf{Case 1:}\label{1st-case-ref} In this situation, between the
jump times the process $\overline{\Xa}_{t}$ evolves as the process
$\Xa_{t}$. The rate of the jumps is given by the function $\Ua_{t}$.
In other words, the jump times $(T_{n})_{n\geq0}$ are given by the
following recursive formulae 
\[
T_{n+1}=\inf{\left\{ t\geq T_{n}~:~\int_{T_{n}}^{t}\Ua_{s}(\overline{\Xa}_{s})~ds\geq e_{n}\right\} }
\]
where $T_{0}=0$, and $(e_{n})_{n\geq0}$ stands for a sequence of
i.i.d. exponential random variables with unit parameter. At the jump
time $T_{n}$ the process $\overline{\Xa}_{T_{n}-}=x$ jumps to new
site $\overline{\Xa}_{T_{n}}=y$ randomly chosen with the distribution
$\mu_{T_{n}-}(dy)$.

For any $f\in D(L)$ we also have that 
\[
\Gamma_{\widehat{L}_{t,\mu_{t}}}(f,f)(x)=\widehat{L}_{t,\mu_{t}}\left((f-f(x))^{2}\right)(x)=\Ua_{t}(x)~\int~[f(y)-f(x)]^{2}~\mu_{t}(dy).
\]

In this situation, an explicit expression of the time inhomogeneous
semigroup $\overline{\Pa}_{s,t,\mu_{s}}$, $s\leq t$, of the process
$\overline{\Xa}_{t}$ is provided by the following formula 
\begin{eqnarray*}
\overline{\Pa}_{s,t,\mu_{s}}(f)(x) & = & \EE\left(f(\overline{X}_{t})~\left|~\overline{X}_{s}=x\right.\right)\\
 & = & Q_{s,t}(1)(x)~\Phi_{s,t}(\delta_{x})(f)+\left(1-Q_{s,t}(1)(x)\right)~\Phi_{s,t}(\mu_{s})(f)\\
 & = & Q_{s,t}(f)(x)+\left(1-Q_{s,t}(1)(x)\right)~\Phi_{s,t}(\mu_{s})(f).
\end{eqnarray*}

We let $\overline{\Pa}_{t_{n},t_{n+1},\mu_{t_{n}}}^{(m)}$ and $\Phi_{t_{n},t_{n+1}}^{(m)}$
be the Markov transition and the transformation of probability measures
defined as $\overline{\Pa}_{s,t,\mu_{s}}$ and $\Phi_{t_{n},t_{n+1}}$
replacing $Q_{t_{n},t_{n+1}}$ by the integral operator 
\begin{eqnarray*}
Q_{t_{n},t_{n+1}}^{(m)}(f)(x) & = & e^{-\Ua_{t_{n}}(x)/m}~\EE\left(f(\Xa_{t_{n+1}})~|~ \Xa_{t_{n}}=x\right).
\end{eqnarray*}
Under the assumptions of theorem~\ref{theo1}, using elementary calculations
we prove that 
\begin{equation}
\overline{\Pa}_{t_{n},t_{n+1},\mu_{t_{n}}}^{(m)}=\overline{\Pa}_{t_{n},t_{n+1},\mu_{t_{n}}}+\frac{1}{m}~\Ra_{t_{n},t_{n+1},\mu_{t_{n}}}^{(m)}\label{sg-m}
\end{equation}
with some remainder signed measures $\Ra_{t_{n},t_{n+1},\mu_{t_{n}}}^{(m)}$
such that $\sup_{m\geq1}{\left\Vert \Ra_{t_{n},t_{n+1},\mu_{t_{n}}}^{(m)}\right\Vert _{\rm tv}}\leq c_{t_{n}}$,
for some finite constant whose values only depend on the potential
function $\Ua_{t}$. 

\item \textbf{Case 2:} In this situation, between jump times the process
$\overline{\Xa}_{t}$ evolves as the process $\Xa_{t}$. The rate
of the jumps is given by the parameter $\mu_{t}(\Va_{t})$. In other
words, the jump times $(T_{n})_{n\geq0}$ are given by the following
recursive formulae 
\[
T_{n+1}=\inf{\left\{ t\geq T_{n}~:~\int_{T_{n}}^{t}\mu_{s}(\Va_{s})~ds\geq e_{n}\right\} }
\]
where $T_{0}=0$, and $(e_{n})_{n\geq0}$ stands for a sequence of
i.i.d. exponential random variables with unit parameter. At the jump
time $T_{n}$ the process $\overline{\Xa}_{T_{n}-}=x$ jumps to new
site $\overline{\Xa}_{T_{n}}=y$ randomly chosen with the distribution
$\Psi_{\Va_{T_{n}-}}(\mu_{T_{n}-})(dy)$.

For any $f\in D(L)$ we also have that 
\[
\Gamma_{\widehat{L}_{t,\mu_{t}}}(f,f)(x)=\widehat{L}_{t,\mu_{t}}\left((f-f(x))^{2}\right)(x)=\int~[f(y)-f(x)]^{2}~\Va_{t}(y)~\mu_{t}(dy).
\]

\item \textbf{Case 3:} In this case, between jump times the process $\overline{\Xa}_{t}$
evolves as the process $\Xa_{t}$. The rate of the jumps is given
by the function 
\begin{eqnarray*}
\Wa_{t,\mu_{t}}(x) & := & \mu_{t}((\Va_{t}-\Va_{t}(x))_{+})\\
 & = & \mu_{t}(\Va_{t}~1_{\Va_{t}\geq\Va_{t}(x)})-\mu_{t}(\Va_{t}\geq\Va_{t}(x))~\Va_{t}(x).
\end{eqnarray*}
In other words, the jump times $(T_{n})_{n\geq0}$ are given by the
following recursive formulae 
\[
T_{n+1}=\inf{\left\{ t\geq T_{n}~:~\int_{T_{n}}^{t}\Wa_{t,\mu_{t}}(\overline{\Xa}_{s})~ds\geq e_{n}\right\} }
\]
where $T_{0}=0$, and $(e_{n})_{n\geq0}$ stands for a sequence of
i.i.d. exponential random variables with unit parameter. At the jump
time $T_{n}$ the process $\overline{\Xa}_{T_{n}-}=x$ jumps to new
site $\overline{\Xa}_{T_{n}}=y$ randomly chosen with the distribution
$\Psi_{(\Va_{T_{n}}-\Va_{T_{n}}(x))_{+}}(\mu_{T_{n}-})$.

For any $f\in D(L)$ we also have that 
\[
\Gamma_{\widehat{L}_{t,\mu_{t}}}(f,f)(x)=\widehat{L}_{t,\mu_{t}}\left((f-f(x))^{2}\right)(x)=\int~[f(y)-f(x)]^{2}~(\Va_{t}(y)-\Va_{t}(x))_{+}~\mu_{t}(dy)
\]
so that 
\begin{equation}
\mu_{t}\left[\Gamma_{\widehat{L}_{t,\mu_{t}}}(f,f)\right]=\int~[f(y)-f(x)]^{2}~(\Va_{t}(y)-\Va_{t}(x))_{+}~\mu_{t}(dx)\mu_{t}(dy).\label{c3}
\end{equation}

\end{itemize}
\textcolor{red}{}

We end this section with another McKean interpretation model combining
cases 1 and 2, as an alternative to the generator described in the
latter case. First, using the fact that 
\[
\mu_{t}\left(\left[\Va_{t}-\mu_{t}(\Va_{t})\right]_{+}-\left[\Va_{t}-\mu_{t}(\Va_{t})\right]_{-}\right)=\mu_{t}\left(\left[\Va_{t}-\mu_{t}(\Va_{t})\right]\right)=0
\]
we prove the following decompositions 
\begin{eqnarray*}
\mu_{t}(\Va_{t}f)-\mu_{t}(\Va_{t})\mu_{t}(f) & = & \mu_{t}\left(\left[\Va_{t}-\mu_{t}(\Va_{t})\right]f\right)\\
 & = & \mu_{t}\left(\left[\Va_{t}-\mu_{t}(\Va_{t})\right]_{+}f\right)-\mu_{t}\left(\left[\Va_{t}-\mu_{t}(\Va_{t})\right]_{-}f\right)\\
 & = & \mu_{t}\left(\left[\Va_{t}-\mu_{t}(\Va_{t})\right]_{+}\left[f-\mu_{t}(f)\right]\right)\\
 &  & \hskip3cm-\mu_{t}\left(\left[\Va_{t}-\mu_{t}(\Va_{t})\right]_{-}\left[f-\mu_{t}(f)\right]\right).
\end{eqnarray*}
Using the same line of arguments as those used in cases 1 and 2, this
implies that 
\[
\mu_{t}(\Va_{t}f)-\mu_{t}(\Va_{t})\mu_{t}(f)=\mu_{t}(\widehat{L}_{t,\mu_{t}}(f))\quad\mbox{{\rm with}}\quad\widehat{L}_{t,\mu_{t}}=\widehat{L}_{t,\mu_{t}}^{+}+\widehat{L}_{t,\mu_{t}}^{-}
\]
where the pair of interacting jump generators is given by 
\[
\widehat{L}_{t,\mu_{t}}^{-}(f)(x)=\left[\Va_{t}(x)-\mu_{t}(\Va_{t})\right]_{-}~\int~[f(y)-f(x)]~\mu_{t}(dy)
\]
and 
\[
\widehat{L}_{t,\mu_{t}}^{+}(f)(x)=~\int~[f(y)-f(x)]~\left[\Va_{t}(y)-\mu_{t}(\Va_{t})\right]_{+}~\mu_{t}(dy).
\]
In this situation, for any $f\in D(L)$ we also have that 
\begin{eqnarray*}
\Gamma_{\widehat{L}_{t,\mu_{t}}}(f,f)(x) & = & \Gamma_{\widehat{L}_{t,\mu_{t}}^{+}}(f,f)(x)+\Gamma_{\widehat{L}_{t,\mu_{t}}^{-}}(f,f)(x)\\
 & = & ~\int~[f(y)-f(x)]^{2}~\left(\left[\Va_{t}(y)-\mu_{t}(\Va_{t})\right]_{+}+\left[\Va_{t}(x)-\mu_{t}(\Va_{t})\right]_{-}\right)~\mu_{t}(dy)
\end{eqnarray*}
so that 
\[
\mu_{t}\left[\Gamma_{\widehat{L}_{t,\mu_{t}}}(f,f)\right]=\int~[f(y)-f(x)]^{2}~\left|\Va_{t}(y)-\mu_{t}(\Va_{t})\right|~\mu_{t}(dx)\mu_{t}(dy).
\]

\subsection{Mean field particle interpretation models}

\label{mean-fieldC} The mean field $N$-particle model $\xi_{t}:=\left(\xi_{t}^{i}\right)_{1\leq i\leq N}$
associated with a given collection of generators $L_{t,\mu_{t}}$
satisfying the weak equation (\ref{WEAK}) is a Markov process in
$E^{N}$ with infinitesimal generator given by the following formulae
\begin{equation}
\La_{t}(F)(x^{1},\ldots,x^{N}):=\sum_{1\leq i\leq N}L_{t,m(x)}^{(i)}(F)(x^{1},\ldots,x^{i},\ldots,x^{N})\quad\mbox{{\rm with}}\quad m(x):=\frac{1}{N}\sum_{1\leq i\leq N}\delta_{x^{i}}\label{defmeanfield}
\end{equation}
for sufficiently regular functions $F$ on $E^{N}$, and for any $x=(x^{i})_{1\leq i\leq N}\in E^{N}$.
In the above formulae, $L_{t,m(x)}^{(i)}$ stands for the operator
$L_{t,m(x)}$ acting on the function $x^{i}\mapsto F(x^{1},\ldots,x^{i},\ldots,x^{N})$.

Before entering into the description of the particle model associated
with the three cases presented in section~\ref{mckean-sec}, we provide
a brief discussion of the convergence analysis of these stochastic
models. Firstly, we recall that 
\[
dF(\xi_{t})=\La_{t}(F)(\xi_{t})~dt+d\Ma_{t}(F)
\]
for some martingale $\Ma_{t}(\varphi)$ with increasing process given
by 
\[
\langle\Ma(F)\rangle_{t}:=\int_{0}^{t}\Gamma_{\La_{s}}\left(F,F\right)(\xi_{s})~ds.
\]
In the above we denote by $\Gamma_{\La_{s}}$ the carré du champ operator
associated with $\La_{s}$, and defined by 
\[
\Gamma_{\La_{s}}\left(F,F\right)(x):=\La_{s}\left[\left(F-F(x)\right)^{2}\right](x)=\La_{s}(F^{2})(x)-F(x)\La_{s}(F)(x).
\]
For empirical test functions of the following form $F(x)=m(x)(f)$,
with $f\in D(L)$, we find that 
\[
\La_{s}(F)(x)=m(x)(L_{s,m(x)}(f))\quad\mbox{{\rm and}}\quad\Gamma_{\La_{s}}\left(\varphi,\varphi\right)(x)=\frac{1}{N}~m(x)\left(\Gamma_{L_{s,m(x)}}(f,f)\right).
\]
From this discussion, if we set $\mu_{t}^{N}=\frac{1}{N}\sum_{1\leq i\leq N}\delta_{\xi_{t}^{i}}$,
then we find that 
\[
d\mu_{t}^{N}(f)=\mu_{t}^{N}(L_{t,\mu_{t}^{N}}(f))~dt+\frac{1}{\sqrt{N}}~dM_{t}^{N}(f)
\]
for any $f\in D(L)$, with the martingale $M_{t}^{N}(f)=\sqrt{N}\Ma_{t}(F)$
with angle bracket given by 
\[
\langle M^{N}(f)\rangle_{t}:=\int_{0}^{t}\mu_{s}^{N}\left(\Gamma_{L_{s,\mu_{s}^{N}}}(f,f)\right)~ds.
\]
A more explicit description of the r.h.s. terms in the above can be
given in the three cases discussed in section~\ref{mckean-sec}.
For instance, in the third case, using formula (\ref{c3}) we find
that 
\[
\langle M^{N}(f)\rangle_{t}:=\int_{0}^{t}\int[f(y)-f(x)]^{2}~(\Va_{s}(y)-\Va_{s}(x))_{+}~\mu_{s}^{N}(dx)\mu_{s}^{N}(dy)~ds.
\]

We conclude that $\mu_{t}^{N}$ ``almost solve'', as $N\uparrow\infty$,
the nonlinear evolution equation (\ref{WEAK}). For a more thorough
discussion of these continuous time models, we refer to the reader
to the review article~\cite{dh-2010}, and the references therein.

By construction, the generator $\La_{t}$ associated with the nonlinear
model (\ref{WEAK}) is decomposed into a mutation generator $\La_{t}^{\rm mut}$
and an interacting jump generator $\La_{t}^{\rm jump}$ 
\[
\La_{t}=\La_{t}^{\rm mut}+\La_{t}^{\rm jump}
\]
with $\La_{t}^{\rm mut}$ and $\La_{t}^{\rm jump}$ defined by
\begin{eqnarray*}
\La_{t}^{\rm mut}(F)(x) & = & \sum_{1\leq i\leq N}L_{t}^{(i)}(F)(x^{1},\ldots,x^{i},\ldots,x^{N})\\
\La_{t}^{\rm jump}(F)(x) & = & \sum_{1\leq i\leq N}\widehat{L}_{t,m(x)}^{(i)}(F)(x^{1},\ldots,x^{i},\ldots,x^{N}).
\end{eqnarray*}
The mutation generator $\La_{t}^{\rm mut}$ describes the evolution
of the particles between the jumps. Between jumps, the particles evolve
independently with $L_{t}$-motions in the sense that they explore
the state space as independent copies of the process $\Xa_{t}$ with
generator $L_{t}$. The jump transition depends on the form of the
generator $\widehat{L}_{t,\mu_{t}}$.
\begin{itemize}
\item \textbf{Case 1:} In this situation the jump generator is given by
\begin{eqnarray*}
\La_{t}^{\rm jump}(F)(x) & = & \sum_{1\leq i\leq N}\Ua_{t}(x^{i})~\int[F(\theta_{u}^{i}(x))-F(x)]~m(x)(du)
\end{eqnarray*}
with the population mappings $\theta_{u}^{i}$ defined below 
\[
\theta_{u}^{i}~:~x\in E^{N}\mapsto\theta_{u}^{i}(x)=(x^{1},\ldots,x^{i-1},\underbrace{u}_{\mbox{ i-th}},x^{i+1},\ldots,x^{N})\in E^{N}.
\]
The quantity $\Ua_{t}(x_{t}^{i})$ represents the jump rate of the
$i$-th particle $\xi_{t}^{i}$. More precisely, if we denote by $T_{n}^{i}$
the $n$-th jump time of $\xi_{t}^{i}$, we have 
\begin{equation}
T_{n+1}^{i}=\inf{\left\{ t\geq T_{n}^{i}~:~\int_{T_{n}^{i}}^{t}\Ua_{s}(\xi_{s}^{i})ds\geq e_{n}^{i}\right\} }\label{TC1}
\end{equation}
where $(e_{n}^{i})_{1\leq i\leq N,n\in\NN}$ stands for a sequence
of i.i.d. exponential random variables with unit parameter. At the
jump time $T_{n}^{i}$ the process $\xi_{T_{n}^{i}-}^{i}=x^{i}$ jumps
to new site $\xi_{T_{n}^{i}}^{i}=u$ randomly chosen with the distribution
$m(\xi_{T_{n}^{i}-})(du)$. In other words, at the jump time the $i$-th
particle jumps to a new state randomly chosen in the current population.

The probabilistic interpretation of the jump generator is not unique.
For instance, it is easily checked that $\La_{t}^{\rm jump}$ can
be rewritten in the following form 
\begin{eqnarray*}
\La_{t}^{\rm jump}(F)(x) & = & \lambda_{t}(x)\int[F(y)-F(x)]~\Pa_{t}(x,dy)
\end{eqnarray*}
with the population jump rate $\lambda_{t}(x)$ and the Markov transition
$\Pa_{t}(x,dy)$ on $E^{N}$ given below 
\[
\lambda_{t}(x):=Nm(x)\left(\Ua_{t}\right)\quad\mbox{{\rm and}}\quad\Pa_{t}(x,dy)=\sum_{1\leq i\leq N}\frac{\Ua_{t}(x^{i})}{\sum_{1\leq i^{\prime}\leq N}\Ua_{t}(x^{i^{\prime}})}~\frac{1}{N}\sum_{1\leq j\leq N}\delta_{\theta_{x^{j}}^{i}(x)}(dy).
\]
In this interpretation, the individual jumps are replaced by population
jumps at rate $\lambda_{t}(\xi_{t})$. More precisely, the jump times
$T_{n}$ of the whole population are defined by 
\[
T_{n+1}=\inf{\left\{ t\geq T_{n}~:~\int_{T_{n}}^{t}\left[\sum_{1\leq i\leq N}\Ua_{s}(\xi_{s}^{i})\right]ds\geq e_{n}\right\} }
\]
where $(e_{n})_{n\in\NN}$ stands for a sequence of i.i.d. exponential
random variables with unit parameter. At the jump time $T_{n}$ the
population $\xi_{T_{n}-}=x$ jumps to new population $\xi_{T_{n}}=y$
randomly chosen with the distribution $\Pa_{T_{n}-}(\xi_{T_{n}-},dy)$.
In other words, at the jump time $T_{n}$, we select randomly a state
$\xi_{T_{n}-}^{i}$ with a probability proportional to $\Ua_{t}(\xi_{T_{n}-}^{i})$,
and we replace this state by a randomly chosen state $\xi_{T_{n}-}^{j}$
in the population, with $1\leq j\leq N$. We end this description
with an alternative interpretation when $\|\Ua_{t}\|\leq C$ for some
finite constant $C<\infty$. In this situation, we clearly have $\|\lambda_{t}\|\leq NC$
and 
\[
\La_{t}^{\rm jump}(F)(x)=\lambda^{\prime}\int[F(y)-F(x)]~\Pa_{t}^{\prime}(x,dy)
\]
with the jump rate $\lambda^{\prime}$ and the Markov jump transitions
$\Pa_{t}^{\prime}$ defined below 
\[
\lambda^{\prime}=NC\quad\mbox{{\rm and}}\quad\Pa_{t}^{\prime}(x,dy):=\frac{\lambda_{t}(x)}{NC}~\Pa_{t}(x,dy)+\left(1-\frac{\lambda_{t}(x)}{NC}\right)~\delta_{x}(dy).
\]
In this interpretation, the population jump times $T_{n}$ arrive
at the higher rate $\lambda^{\prime}=NC$. At the jump time $T_{n}$
the population $\xi_{T_{n}-}=x$ jumps to new population $\xi_{T_{n}}=y$
randomly chosen with the distribution $\Pa_{T_{n}-}^{\prime}(\xi_{T_{n}-},dy)$.

In the models described above, as usual, between the jump times $T_{n}$
of the population every particle evolves independently with $L_{t}$-motions. 

\item \textbf{Case 2 :} In this situation, the jump generator is given by
\begin{eqnarray*}
\La_{t}^{\rm jump}(F)(x) & = & \sum_{1\leq i\leq N}m(x)(\Va_{t})~\int[F(\theta_{u}^{i}(x))-F(x)]~\Psi_{\Va_{t}}\left(m(x)\right)(du).
\end{eqnarray*}
The particles have a common jump rate given by the empirical average
$m(\xi_{t})(\Va_{t})$. In other words, the jump times $T_{n}^{i}$
of a particle $\xi_{t}^{i}$ are given by the following recursive
formulae 
\[
T_{n+1}^{i}=\inf{\left\{ t\geq T_{n}^{i}~:~\int_{T_{n}^{i}}^{t}m(\xi_{s})(\Va_{s})~ds\geq e_{n}^{i}\right\} }
\]
where $(e_{n}^{i})_{1\leq i\leq N,n\geq0}$ stands for a sequence
of i.i.d. exponential random variables with unit parameter. At the
jump time $T_{n}^{i}$ the process $\xi_{T_{n}-}^{i}=x^{i}$ jumps
to new site $\xi_{T_{n}}^{i}=u$ randomly chosen with the weighted
distribution $\Psi_{\Va_{T_{n}-}}(m(\xi_{T_{n}-}))(du)$.

As mentioned in the first case, the probabilistic interpretation of
the jump generator is not unique. In this situation, it is easily
checked that $\La_{t}^{\rm jump}$ can be rewritten in the following
form 
\begin{eqnarray*}
\La_{t}^{\rm jump}(F)(x) & = & \lambda_{t}(x)~\int[F(y)-F(x)]~\Pa_{t}(x,dy)
\end{eqnarray*}
with the population jump rate $\lambda_{t}$ and the Markov transition
$\Pa_{t}(x,dy)$ on $E^{N}$ defined below 
\[
\lambda_{t}(x):=Nm(x)(\Va_{t})\quad\mbox{{\rm and}}\quad\Pa_{t}(x,dy)=\frac{1}{N}\sum_{1\leq i\leq N}\sum_{1\leq j\leq N}\frac{\Va_{t}(x^{j})}{\sum_{1\leq j^{\prime}\leq N}\Va_{t}(x^{j\prime})}\delta_{\theta_{x^{j}}^{i}(x)}(dy).
\]
The description of the evolution of the population model follows the
same lines as the ones given in case 1. 

\item \textbf{Case 3 :} In this situation, the jump generator is given by
\begin{eqnarray*}
\La_{t}^{\rm jump}(F)(x) & = & \sum_{1\leq i\leq N}\int[F(\theta_{u}^{i}(x))-F(x)]~(\Va_{t}(u)-V_{t}(x^{i}))_{+}~m(x)(du)\\
 & = & \sum_{1\leq i\leq N}m(x)((\Va_{t}-\Va_{t}(x^{i}))_{+})~\int[F(\theta_{u}^{i}(x))-F(x)]~\Psi_{(\Va_{t}-V_{t}(x^{i}))_{+}}(m(x))(du).
\end{eqnarray*}

\end{itemize}
In this interpretation, the jump rate of the $i$-th particle is given
by the average potential variation of the particle with higher values
\[
m(x)((\Va_{t}-\Va_{t}(x^{i}))_{+})=\frac{1}{N}\sum_{1\leq j\leq N}~1_{\left\{ \Va_{t}(x^{j})>\Va_{t}(x^{i})\right\} }~\left(\Va_{t}(x^{j})-\Va_{t}(x^{i})\right).
\]
More precisely, if we denote by $T_{n}^{i}$ the $n$-th jump time
of $\xi_{t}^{i}$, we have 
\[
T_{n+1}^{i}=\inf{\left\{ t\geq T_{n}^{i}~:~\int_{T_{n}^{i}}^{t}m(\xi_{s})((\Va_{s}-\Va_{s}(\xi_{s}^{i}))_{+})ds\geq e_{n}^{i}\right\} }
\]
where $(e_{n}^{i})_{1\leq i\leq N,n\in\NN}$ stands for a sequence
of i.i.d. exponential random variables with unit parameter. At the
jump time $T_{n}^{i}$ the particle $\xi_{T_{n}^{i}-}^{i}=x^{i}$
jumps to new site $\xi_{T_{n}^{i}}^{i}=u$ randomly chosen with the
distribution 
\[
\Psi_{(\Va_{T_{n}^{i}-}-\Va_{T_{n}^{i}-}(x^{i}))_{+}}(m(\xi_{T_{n}^{i}-}))(du)\propto\sum_{1\leq j\leq N}~1_{\left\{ \Va_{T_{n}^{i}-}(\xi_{T_{n}^{i}-}^{j})>\Va_{T_{n}^{i}-}(x^{i})\right\} }~\left(\Va_{T_{n}^{i}-}(\xi_{T_{n}^{i}-}^{j})-\Va_{T_{n}^{i}-}(x^{i})\right)~\delta_{\xi_{T_{n}^{i}-}^{j}}(du).
\]
In other words, we choose randomly a new site $\xi_{T_{n}^{i}}^{i}=\xi_{T_{n}^{i}-}^{j}$,
among the ones with higher potential value with a probability proportional
to the difference of potential $\left(\Va_{T_{n}^{i}-}(\xi_{T_{n}^{i}-}^{j})-\Va_{T_{n}^{i}-}(\xi_{T_{n}^{i}-}^{i})\right)$.

As the first two cases discussed above, we can also interpret this
jump generator at the level of the population. In this interpretation
we have 
\[
\La_{t}^{\rm jump}(F)(x)=\lambda_{t}(x)~\int[F(y)-F(x)]~\Pa_{t}(x,dy)
\]
with the population jump rate 
\[
\lambda_{t}(x)=N\int m(x)(du)~m(x)(dv)~(\Va_{t}(u)-\Va_{t}(v))_{+}
\]
and the population jump transition 
\begin{eqnarray*}
\Pa_{t}(x,dy) & = & \sum_{1\leq i,j\leq N}\frac{(\Va_{t}(x^{j})-V_{t}(x^{i}))_{+}}{\sum_{1\leq i^{\prime},j^{\prime}\leq N}(\Va_{t}(x^{j^{\prime}})-\Va_{t}(x^{i^{\prime}}))_{+}}~~\delta_{\theta_{x^{j}}^{i}(x)}(dy)\ .
\end{eqnarray*}

\begin{rem} In case \textbf{(D)}, the reference Markov process $\Xa_{t}=X_{\lfloor t\rfloor}$
has deterministic and fixed time jumps on integer times so that the
generator approach developed above does not apply directly. Nevertheless
their probabilistic interpretation is defined in the same way:

Between the jumps, the process $\overline{X}_{t}$ evolves as $\Xa_{t}$,
and the $N$ particles explore the state space as independent copies
of the process $\Xa_{t}$. The rate of the jumps and their random
spatial location are defined using the same interpretations as the
ones given above.

The stochastic modeling and the analysis of these continuous time
models and their particle interpretations can be developed using the
semigroup techniques provided in~\cite{dm2007}. \end{rem}

\section{Discrete time models}

\label{sec-disc-mod}

\subsection{McKean models and Feynman-Kac semigroups}

\label{mckean-secD}

As in the continuous time case, these discrete time evolution equations
(\ref{WEAKD}) can be interpreted as the evolution of the laws defined
by $\mbox{{\rm Law}}\left(\overline{X}_{t_{n}}\right)=\mu_{t_{n}}^{(m)}$
of a time inhomogeneous Markov process $\overline{X}_{t_{n}}$ with
Markov transitions generators $\Ka_{t_{n},t_{n+1},\mu_{t_{n}}^{(m)}}$
that depend on the distribution of the random states at the previous sub-integer mesh time increment. This probabilistic
model is also called the McKean interpretation of the evolution equation
(\ref{WEAKD}) in terms of a time inhomogeneous Markov chain. By construction,
the elementary transitions of the Markov chain $\overline{X}_{t_{n}}\leadsto\overline{X}_{t_{n+1}}$
are decomposed into two separate transitions $\overline{X}_{t_{n}}\leadsto\widehat{X}_{t_{n}}\leadsto\overline{X}_{t_{n+1}}$.

First, the state $\overline{X}_{t_{n}}=x$ jumps to a new location
$\widehat{X}_{t_{n}}=y$ randomly chosen with the Markov transition
$\Sa_{t_{n},\mu_{t_{n}}^{(m)}}(x,dy)$, given by one of the three cases presented in section 2.2 when applied to the particle approximation of the measure $\mu_{t_{n}}^{(m)}$ at mesh time increment $t_{n}$. Then, the selected state $\widehat{X}_{t_{n}}=y$
evolves to a new site $\overline{X}_{t_{n+1}}=z$ according to the
Markov transition $\Ma_{t_{n},t_{n+1}}(y,dz)$.

Next, we recall some basic properties of the semigroup $\Phi_{t_{p},t_{n}}^{(m)}$
of the flow of measures $\mu_{t_{n}}^{(m)}$. By construction, we
have 
\[
\Phi_{t_{p},t_{n}}^{(m)}(\mu_{t_{p}}^{(m)})(f)={\mu_{t_{p}}^{(m)}Q_{t_{p},t_{n}}^{(m)}(f)}/{\mu_{t_{p}}^{(m)}Q_{t_{p},t_{n}}^{(m)}(1)}
\]
with Feynman-Kac semigroup $Q_{t_{p},t_{n}}^{(m)}$ defined by 
\begin{eqnarray}
Q_{t_{p},t_{n}}^{(m)}(f)(x) & = & \EE\left(f(\Xa_{t_{n}})~\prod_{p\leq q<n}e^{\Va_{t_{q}}(\Xa_{t_{q}})/m}~\left|~\Xa_{t_{p}}=x\right.\right).\label{Q-ref}
\end{eqnarray}

We notice that when we consider the Boltzmann-Gibbs transformation associated with the potential
function $G_{t_{p},t_{n}}^{(m)}=Q_{t_{p},t_{n}}^{(m)}(1)$ then the semigroup of the flow of measures $\mu_{t_{n}}^{(m)}$ can be expressed according to 
\begin{equation}
\Phi_{t_{p},t_{n}}^{(m)}(\mu_{t_{p}}^{(m)})=\Psi_{G_{t_{p},t_{n}}^{(m)}}(\mu_{t_{p}}^{(m)})P_{t_{p},t_{n}}^{(m)}\quad\mbox{{\rm with}}\quad P_{t_{p},t_{n}}^{(m)}(f)={Q_{t_{p},t_{n}}^{(m)}(f)}/{Q_{t_{p},t_{n}}^{(m)}(1)}\label{dec-ref1}.
\end{equation}

\begin{defi} We consider the integral operators 
\[
L_{t_{n},\mu}^{(m)}:=\Ka_{t_{n},t_{n+1},\mu}-Id,\qquad L_{t_{n}}^{(m)}:=\Ma_{t_{n},t_{n+1}}-Id\quad\mbox{and}\quad\widehat{L}_{t_{n},\mu}^{(m)}:=\Sa_{t_{n},\mu}-Id.
\]
\end{defi}

\begin{lem}\label{prop-dec} We have the decomposition 
\begin{eqnarray*}
L_{t_{n},\mu}^{(m)} & = & L_{t_{n}}^{(m)}+\widehat{L}_{t_{n},\mu}^{(m)}+\widehat{L}_{t_{n},\mu}^{(m)}L_{t_{n}}^{(m)}.
\end{eqnarray*}
In addition, for any $f\in\Ba_{b}(E)$ $\mu\in\Pa(E)$, and any $x\in E$,
we have 
\[
\Ka_{t_{n},t_{n+1},\mu}\left(\left[f-\Ka_{t_{n},t_{n+1},\mu}(f)(x)\right]^{2}\right)(x)=\Gamma_{L_{t_{n},\mu}^{(m)}}(f,f)(x)-\left(L_{t_{n},\mu}^{(m)}(f)(x)\right)^{2}.
\]
\end{lem} \proof Using the decomposition 
\begin{eqnarray*}
L_{t_{n},\mu}^{(m)} & = & \left(\Ma_{t_{n},t_{n+1}}-Id\right)+\left(\Sa_{t_{n},\mu}-Id\right)+\left(\Sa_{t_{n},\mu}-Id\right)\left(\Ma_{t_{n},t_{n+1}}-Id\right)
\end{eqnarray*}
we readily check the first assertion. We prove the second decomposition
using the fact that 
\[
\left[\Ka_{t_{n},t_{n+1},\mu}(f)\right]^{2}=\left(L_{t_{n},\mu}^{(m)}(f)\right)^{2}+f^{2}-2fL_{t_{n},\mu}^{(m)}(f).
\]
This ends the proof of the lemma. \cqfd

In the further development of this section $c_{t_{n}}<\infty$ stands
for some generic finite constant whose values may vary from line to
line. For any function $f\in D(L)$, such that $L_{t}(f)\in\Ca^{1}([t_{n},t_{n+1}],D(L))$,
we also define 
\begin{equation}
\|f\|_{t_{n}}:=\|f\|+\sup_{t_{n}\leq t\leq t_{n+1}}\left(\left\Vert {\partial L_{t}(f)}/{\partial t}\right\Vert +\left\Vert L_{t}(f)\right\Vert +\left\Vert L_{t}^{2}(f)\right\Vert \right).\label{normftn}
\end{equation}

\begin{prop}\label{lem-ref-D} In case \textbf{(D)}, we have 
\[
L_{t_{n}}^{(m)}=1_{\NN}(t_{n+1})~\left(M_{t_{n+1}}-Id\right)\quad\mbox{ and}\quad L_{t_{n},\mu}^{(m)}=\widehat{L}_{t_{n},\mu}^{(m)}+1_{\NN}(t_{n+1})~\Sa_{t_{n},\mu}\left(M_{t_{n+1}}-Id\right).
\]

In case \textbf{(C)}, we have the first order expansion 
\begin{equation}
L_{t_{n}}^{(m)}(f)=L_{t_{n}}(f)~ \frac{1}{m}+R_{t_{n}}(f)~\frac{1}{m^{2}}.\label{1st-dec-diff}
\end{equation}
for any function $f\in D(L)$, such that $L_{t}(f)\in\Ca^{1}([t_{n},t_{n+1}],D(L))$,
with some remainder operator $R_{t_{n}}$ such that $\|R_{t_{n}}(f)\|\leq c_{t_{n}}~\|f\|_{t_{n}}$.
Furthermore, we have the first order expansion 
\begin{equation}
L_{t_{n},\mu}^{(m)}(f)=\frac{1}{m}~L_{t_{n},\mu}(f)+\frac{1}{m^{2}}~R_{t_{n},\mu}(f)\label{1st-dec-diff-C}
\end{equation}
with some second order remainder term $R_{t_{n},\mu}(f)$ such that
$\sup_{\mu\in\Pa(E)}{\left\Vert R_{t_{n},\mu}(f)\right\Vert }\leq c_{t_{n}}\|f\|_{t_{n}}$.
In addition, we have 
\begin{equation}
\Ka_{t_{n},t_{n+1},\mu}\left(\left[f-\Ka_{t_{n},t_{n+1},\mu}(f)(x)\right]^{2}\right)(x)=\Gamma_{L_{t_{n},\mu}}(f,f)(x)~ \frac{1}{m}+\Ra_{t_{n},\mu}\left(f,f\right)(x) ~\frac{1}{m^{2}}\label{ref-Ka}
\end{equation}
with some remainder operator s.t. $\sup_{\mu\in\Pa(E)}{\left\Vert \Ra_{t_{n},\mu}\left(f,f\right)\right\Vert }\leq c_{t_{n}}~\|f^{2}\|_{t_{n}}$.
\end{prop} The proof of proposition~\ref{lem-ref-D} is provided
in the appendix, on page~\pageref{proof-lem-DC}.

The first order expansions stated in the proposition~\ref{lem-ref-D}
can be used to develop a stochastic perturbation approach to estimate
the deviations of the measures $\mu_{t_{n}}^{(m)}$ around their limiting
values $\mu_{t_{n}}$. Next, we provide an alternative approach based
on the explicit representation (\ref{sg-m}) of the time inhomogeneous
transition of the limiting process $\overline{\Xa}_{t_{n}}$ on the
time mesh sequence $t_{n}$. In the first case discussed on page~\pageref{1st-case-ref},
we have 
\begin{equation}\label{Phi-P-ref}
\mu_{t_{n+1}}^{(m)}=\Phi_{t_{n},t_{n+1}}^{(m)}\left(\mu_{t_{n}}^{(m)}\right):=\Psi_{e^{-\Ua_{t_{n}}/m}}\left(\mu_{t_{n}}^{(m)}\right)\Ma_{t_{n},t_{n+1}}=\mu_{t_{n}}^{(m)}\overline{\Pa}_{t_{n},t_{n+1},\mu_{t_{n}}^{(m)}}^{(m)}
\end{equation}
with the Markov transition 
\begin{equation}\label{Phi-P-ref-M}
\begin{array}{l}
\overline{\Pa}_{t_{n},t_{n+1},\mu_{t_{n}}^{(m)}}^{(m)}(x,dy)\\
\\
=e^{-\Ua_{t_{n}}(x)/m}~\Ma_{t_{n},t_{n+1}}(x,dy)+\left(1-e^{-\Ua_{t_{n}}(x)/m}\right)~\Psi_{e^{-\Ua_{t_{n}}/m}}\left(\mu_{t_{n}}^{(m)}\right)\Ma_{t_{n},t_{n+1}}(dy).
\end{array}
\end{equation}
Using (\ref{sg-m}) we readily find that 
\[
\mu_{t_{n+1}}^{(m)}=\mu_{t_{n}}^{(m)}\overline{\Pa}_{t_{n},t_{n+1},\mu_{t_{n}}^{(m)}}+\frac{1}{m}~\Wa_{t_{n},t_{n+1}}^{(m)}=\Phi_{t_{n},t_{n+1}}\left(\mu_{t_{n}}^{(m)}\right)+\frac{1}{m}~\Wa_{t_{n},t_{n+1}}^{(m)}
\]
with the signed measure 
\[
\Wa_{t_{n},t_{n+1}}^{(m)}:=\mu_{t_{n}}^{(m)}\Ra_{t_{n},t_{n+1},\mu_{t_{n}}}^{(m)}\quad\mbox{{\rm s.t.}}\quad\sup_{m\geq1}{\left\Vert \Wa_{t_{n},t_{n+1}}^{(m)}\right\Vert _{\rm tv}}\leq c_{t_{n}}
\]
for some finite constant whose values only depend on the potential
function $\Ua_{t}$. In summary, we have proven the following first
order local perturbation decompositions 
\begin{eqnarray*}
\mu_{t_{n+1}}^{(m)} & = & \Phi_{t_{n},t_{n+1}}\left(\mu_{t_{n}}^{(m)}\right)+\frac{1}{m}~\Wa_{t_{n},t_{n+1}}^{(m)}\\
\mu_{t_{n+1}} & = & \Phi_{t_{n},t_{n+1}}\left(\mu_{t_{n}}\right)
\end{eqnarray*}
These local expansions allow the use of perturbation theory developed
in section 7.1 of \cite{dm2004} to derive several qualitative estimates
between $\mu_{t_{n}}^{(m)}$ and $\mu_{t_{n}}$ in terms of the stability
properties of the Feynman-Kac semigroup $\Phi_{t_{n},t_{n+1}}$.

\subsection{Mean field particle interpretation models}

\label{mean-fieldD} If we set $\mu_{t_{n}}^{N}=\frac{1}{N}\sum_{1\leq i\leq N}\delta_{\xi_{t_{n}}^{i}}$,
then we have the decomposition 
\[
\mu_{t_{n+1}}^{N}=\mu_{t_{n}}^{N}\Ka_{t_{n},t_{n+1},\mu_{t_{n}}^{N}}+\frac{1}{\sqrt{N}}~W_{t_{n},t_{n+1}}^{N}
\]
with the sequence of empirical random fields $W_{t_{n},t_{n+1}}^{N}$
such that 
\[
\EE\left(W_{t_{n},t_{n+1}}^{N}(f)~|~\xi_{t_{n}}\right)=0
\]
and 
\[
\EE\left(W_{t_{n},t_{n+1}}^{N}(f)^{2}~|~\xi_{t_{n}}\right)=\int\mu_{t_{n}}^{N}(du)~\Ka_{t_{n},t_{n+1},\mu_{t_{n}}^{N}}(u,dv)~\left(f(v)-\Ka_{t_{n},t_{n+1},\mu_{t_{n}}^{N}}(f)(u)\right)^{2}.
\]

As for the continuous time models, we conclude that $\mu_{t_{n}}^{N}$
``almost solve'', as $N\uparrow\infty$, the nonlinear evolution
equation (\ref{WEAKD}). For a more thorough discussion of these local
sampling random field models, we refer the reader to~\cite{dm2004,dh-2010,dhw2012},
and references therein.

By construction, the elementary transitions of the Markov chain $\xi_{t_{n}}\leadsto\xi_{t_{n+1}}$
are decomposed into two separate transitions: 
\begin{equation}
\xi_{t_{n}}\leadsto\widehat{\xi}_{t_{n}}=\left(\widehat{\xi}_{t_{n}}^{i}\right)_{1\leq i\leq N}\leadsto\xi_{t_{n+1}}.\label{TWO}
\end{equation}

First, every particle $\xi_{t_{n}}^{i}=x^{i}$ jumps independently
to a new location $\widehat{\xi}_{t_{n}}^{i}=y^{i}$ randomly chosen
with the Markov transition $\Sa_{t_{n},m(\xi_{t_{n}})}(x^{i},dy^{i})$,
with $1\leq i\leq N$. Following this, each particle $\widehat{\xi}_{t_{n}}^{i}=y^{i}$
evolves independently to a new site $\xi_{t_{n+1}}^{i}=z^{i}$ according
to the Markov transition $\Ma_{t_{n},t_{n+1}}(y^{i},dz^{i})$, with
$1\leq i\leq N$.

In other words, the mutation transition describes the evolution of
the particles between the jumps. Between the jumps, the particles
evolve independently with $\Ma_{t_{n},t_{n+1}}$-motions in the sense
that they explore the state space as independent copies of the process
$\Xa_{t_{n}}$ with Markov transition $\Ma_{t_{n},t_{n+1}}$. The
jump transition can also be interpreted as an acceptance-rejection
transition equipped with a recycling mechanism. In this interpretation,
the mutation transition can be interpreted as a proposal transition.
Notice that the selection type transition is dictated by the choice
of the transition $\Sa_{t_{n},\mu_{t_{n}}}$.

We illustrate these jump type transitions in the first case presented
on page \pageref{case1-ref}. In this situation, we recall that the
selection transition of the $i$-th particle $\xi_{t_{n}}^{i}\leadsto\widehat{\xi}_{t_{n}}^{i}$
is given by the following distribution 
\begin{equation}\label{def-S-ref}
\Sa_{t_{n},\mu_{t_{n}}^{N}}(\xi_{t_{n}}^{i},dy):=e^{-\Ua_{t_{n}}(\xi_{t_{n}}^{i})/m}~\delta_{\xi_{t_{n}}^{i}}(dy)+\left(1-e^{-\Ua_{t_{n}}(\xi_{t_{n}}^{i})/m}\right)~\Psi_{e^{-\Ua_{t_{n}}/m}}(\mu_{t_{n}}^{N})(dy).
\end{equation}

Next, we provide an interpretation of this transition as an acceptance-rejection
scheme with a recycling mechanism. We let $\widetilde{\xi}_{t_{n}}=\left(\widetilde{\xi}_{t_{n}}^{i}\right)_{1\leq i\leq N}$
be a sequence of conditionally independent random variables with common
law 
\[
\Psi_{e^{-\Ua_{t_{n}}/m}}(\mu_{t_{n}}^{N})=\sum_{1\leq i\leq N}\frac{e^{-\Ua_{t_{n}}(\xi_{t_{n}}^{i})/m}}{\sum_{1\leq j\leq N}e^{-\Ua_{t_{n}}(\xi_{t_{n}}^{j})/m}}~\delta_{\xi_{t_{n}}^{i}}.
\]
We also consider a sequence of conditionally independent Bernoulli
random variables with distribution 
\[
\PP\left(\epsilon_{t_{n}}^{i}=1~|~\xi_{t_{n}}\right)=1-\PP\left(\epsilon_{t_{n}}^{i}=0~|~\xi_{t_{n}}\right)=e^{-\Ua_{t_{n}}(\xi_{t_{n}}^{i})/m}.
\]
In this notation, we have that 
\[
\widehat{\xi}_{t_{n}}^{i}=\epsilon_{t_{n}}^{i}~\xi_{t_{n}}^{i}+\left(1-\epsilon_{t_{n}}^{i}\right)~\widetilde{\xi}_{t_{n}}^{i}.
\]
In in other words, the particle $\xi_{t_{n}}^{i}$ is accepted when
$\epsilon_{t_{n}}^{i}=1$; otherwise, it is rejected and replaced
by a particle $\widetilde{\xi}_{t_{n}}^{i}$ randomly chosen with
the updated weighted distribution $\Psi_{e^{-\Ua_{t_{n}}/m}}(\mu_{t_{n}}^{N})$.
The pool of particles that have been accepted from the start provide
a sequence of exact samples. More precisely, it can be easily shown
that 
\[
\mbox{{\rm Law}}(~\xi_{t_{n}}^{i}~|~\forall p:0\leq p<n,~~\epsilon_{t_{p}}^{i}=1)=\mu_{t_{n}}.
\]
See for instance section 1.5.1 in~\cite{dm2004}.

In connection with (\ref{TC1}), we notice that the rejection times
$T_{n}^{i}$ on the time mesh $(t_{q})_{q\geq0}$ can be defined as
follows 
\begin{eqnarray*}
T_{n+1}^{i} & = & \inf{\left\{ t_{p}>T_{n}^{i}~:~\sum_{T_{n}^{i}\leq k\leq t_{p}}\Ua_{t_{k}}(\xi_{t_{k}}^{i})/m\geq e_{n}^{i}\right\} }\\
 & = & \inf{\left\{ t_{p}>T_{n}^{i}~:~\prod_{T_{n}^{i}\leq k\leq t_{p}}e^{-\Ua_{t_{k}}(\xi_{t_{k}}^{i})/m}\leq u_{n}^{i}\right\} }
\end{eqnarray*}
where $(u_{n}^{i})_{1\leq i\leq N,n\in\NN}$ stands for a sequence
of i.i.d.uniform random variables on $]0,1]$, and $e_{n}^{i}:=-\log{u_{n}^{i}}$,
with $1\leq i\leq N,n\in\NN$, is the corresponding sequence of i.i.d.
exponential random variables with unit parameter. We check this claim
using the following observations 
\[
\begin{array}{l}
\PP\left(T_{n+1}^{i}=t_{p}~|~T_{n}^{i},~~\forall~T_{n}^{i}\leq k\leq t_{p}~~\xi_{t_{k}}^{i}\right)\\
\\
=\PP\left(\forall T_{n}^{i}\leq k<t_{p}~~\epsilon_{t_{k}}=1,~~\epsilon_{t_{p}}^{i}=0~|~T_{n}^{i},~~ \forall~T_{n}^{i}\leq k\leq t_{p}~~\xi_{t_{k}}^{i}\right)\\
\\
=\left(\prod_{T_{n}^{i}\leq k<t_{p}}e^{-\Ua_{t_{k}}(\xi_{t_{k}}^{i})/m}\right)\times\left(1-e^{-\Ua_{t_{p}}(\xi_{t_{p}}^{i})/m}\right)\\
\\
=e^{-\sum_{T_{n}^{i}\leq k<t_{p}}\Ua_{t_{k}}(\xi_{t_{k}}^{i})/m}\times\left(1-e^{-\Ua_{t_{p}}(\xi_{t_{p}}^{i})/m}\right)\\
\\
=\PP\left(\sum_{T_{n}^{i}\leq k<t_{p}}\Ua_{t_{k}}(\xi_{t_{k}}^{i})/m<e_{n}^{i}\leq\sum_{T_{n}^{i}\leq k\leq t_{p}}\Ua_{t_{k}}(\xi_{t_{k}}^{i})/m~|~T_{n}^{i},~~\forall~T_{n}^{i}\leq k\leq t_{p}~~\xi_{t_{k}}^{i}\right).
\end{array}
\]
At the jump time $T_{n}^{i}$ the process $\xi_{T_{n}^{i}}^{i}=x^{i}$
jumps to new site $\widehat{\xi}_{T_{n}^{i}}^{i}=u$ randomly chosen
with the distribution 
\begin{equation}\label{recycling-trans}
\Psi_{e^{-\Ua_{T_{n}^{i}}/m}}(\mu_{T_{n}^{i}}^{N}).
\end{equation}

\subsection{An mean field model with uniform recycling}\label{remark:parallelcomputation}
When $m$ is large enough, the recycling distribution (\ref{recycling-trans}) in the Markov transition (\ref{def-S-ref})
is almost equal to $\mu_{T_{n}^{i}}^{N}$. For instance, we have the
total variation estimate 
\[
\left\Vert \Psi_{e^{-\Ua_{T_{n}^{i}}/m}}(\mu_{T_{n}^{i}}^{N})-\mu_{T_{n}^{i}}^{N}\right\Vert _{\rm tv}\leq\left\Vert 1-e^{-\Ua_{T_{n}^{i}}/m}\right\Vert \leq\|\Ua_{T_{n}^{i}}\|/m.
\]
We prove these inequalities using the decomposition 
\[
\Psi_{e^{-\Ua/m}}(\mu)(f)-\mu(f)=\mu\left(\left[1-e^{-\Ua/m}\right]\left[\Psi_{e^{-\Ua/m}}(\mu)(f)-f\right]\right)
\]
which is valid for any bounded functions $\Ua$ and $f$. Hence, to
save computational time we can replace the recycling weighted measure
$\Psi_{e^{-\Ua_{T_{n}^{i}}/m}}(\mu_{T_{n}^{i}}^{N})$ by $\mu_{T_{n}^{i}}^{N}$.
In this situation, the selection transition takes the following form
\[
\widehat{\xi}_{t_{n}}^{i}=\epsilon_{t_{n}}^{i}~\xi_{t_{n}}^{i}+\left(1-\epsilon_{t_{n}}^{i}\right)~\xi_{t_{n}}^{1+\lfloor N\tau_{t_{n}}^{i}\rfloor}
\]
where $\tau_{t_{n}}^{i}$ stands for a sequence of i.i.d. uniform
random variables on $]0,1]$. 

Next, we detail some analysis of the differences between the McKean models with recycling Boltzmann-Gibbs transitions, and the models discussed which utilise uniform recycling.

We let $\widetilde{\Sa}_{t_{n},\mu}$ the collection of selection transitions
defined as in (\ref{def-S-ref}), by replacing $\Psi_{e^{-\Ua_{t_{n}}/m}}(\mu)$ by the measure $\mu$. Furthermore, we denote by $\Phi_{t_p,t_n}^{(m)}$, resp. $\widetilde{\Phi}^{(m)}_{t_p,t_n}$, with $p\leq n$,
the semigroups associated with these flows
$$
\widetilde{\Phi}^{(m)}_{t_p,t_n}\left(\widetilde{\mu}_{t_{p}}^{(m)}\right)=\widetilde{\mu}_{t_{n}}^{(m)}\quad
\mbox{\rm and}\quad
\Phi^{(m)}_{t_p,t_n}\left(\mu_{t_{p}}^{(m)}\right)=\mu_{t_{n}}^{(m)}.
$$
Then, using the decomposition
$$
\left[\Sa_{t_{n},\mu}-\widetilde{\Sa}_{t_{n},\mu}\right](x,dy)=\left(1-e^{-\Ua_{t_{n}}(x)/m}\right)~\left(\Psi_{e^{-\Ua/m}}(\mu)-\mu\right)(dy)
$$
we find that
$$
\sup_{x\in E}{\left\Vert 
\Sa_{t_{n},\mu}(x,\point)-\widetilde{\Sa}_{t_{n},\mu}(x,\point)
\right\Vert _{\rm tv}}\leq \|\Ua_{t_{n}}\|^2/m^2.
$$
Replacing  $\Psi_{e^{-\Ua_{t_{n}}/m}}(\mu)$ by the measure $\mu$ in (\ref{Phi-P-ref-M}),
the evolution equation (\ref{Phi-P-ref}) takes the following form
$$
\widetilde{\mu}_{t_{n+1}}^{(m)}=\widetilde{\Phi}_{t_{n},t_{n+1}}^{(m)}\left(\widetilde{\mu}_{t_{n}}^{(m)}\right):=\widetilde{\mu}_{t_{n}}^{(m)}\widetilde{\Pa}_{t_{n},t_{n+1},\widetilde{\mu}_{t_{n}}^{(m)}}^{(m)}
\quad\mbox{\rm
with}\quad
\widetilde{\Pa}_{t_{n},t_{n+1},\widetilde{\mu}_{t_{n}}^{(m)}}^{(m)}=\widetilde{\Sa}_{t_{n},\widetilde{\mu}_{t_{n}}^{(m)}}\Ma_{t_{n},t_{n+1}}.
$$
From previous estimates, we find that
$$
\widetilde{\Phi}_{t_{n},t_{n+1}}^{(m)}\left(\mu\right)=\Phi_{t_{n},t_{n+1}}^{(m)}\left(\mu\right)+\frac{1}{m^2}~\Ra_{t_{n},t_{n+1}}^{(m)}\left(\mu\right)
$$
with some measures $~\Ra_{t_{n},t_{n+1}}^{(m)}\left(\mu\right)$ s.t. $$
\sup_{\mu\in \Pa(E)}{\left\Vert
\Ra_{t_{n},t_{n+1}}^{(m)}\left(\mu\right)
 \right\Vert _{\rm tv}}\leq \|\Ua_{t_{n}}\|^2~.$$
We end this section with an estimate of the difference between the flow of measures $\widetilde{\mu}_{t_{n+1}}^{(m)}$, and $\mu_{t_{n+1}}^{(m)}$.  We are now in position to state, and to prove the following theorem.
\begin{theo}\label{theo-interm}
For any $n\geq 0$, we have
$$
\left\Vert 
\widetilde{\mu}_{t_{n}}^{(m)}-\mu_{t_{n}}^{(m)}
\right\Vert _{\rm tv}\leq c_{t_n}~/m
$$
for some finite constant $c_{t_n}<\infty$, whose values don't depend on the parameter $m$. 
\end{theo}
\proof
We use the stochastic perturbation analysis developed in section 6.3 in~\cite{dhw2012} (see also chapter 7 in ~\cite{dm2004}).
We denote by $\Phi_{t_p,t_n}^{(m)}$, resp. $\widetilde{\Phi}^{(m)}_{t_p,t_n}$,
with $p\leq n$,
the semigroups associated with these flows
$$
\widetilde{\Phi}^{(m)}_{t_p,t_n}\left(\widetilde{\mu}_{t_{p}}^{(m)}\right)=\widetilde{\mu}_{t_{n}}^{(m)}\quad
\mbox{\rm and}\quad
\Phi^{(m)}_{t_p,t_n}\left(\mu_{t_{p}}^{(m)}\right)=\mu_{t_{n}}^{(m)}
$$
Using the interpolating sequence of measures
$$
0\leq p\leq n\mapsto \Phi_{t_p,t_n}^{(m)}\left(\widetilde{\Phi}^{(m)}_{t_0,t_p}(\widetilde{\mu}_{t_{0}}^{(m)})\right)=\Phi^{(m)}_{t_p,t_n}\left(\widetilde{\mu}_{t_{p}}^{(m)}\right)
$$
from the distribution
$$
\Phi_{t_0,t_n}^{(m)}\left(\widetilde{\Phi}^{(m)}_{t_0,t_0}(\widetilde{\mu}_{t_{0}}^{(m)}))\right)=\Phi^{(m)}_{t_0,t_n}\left(\widetilde{\mu}_{t_{0}}^{(m)})\right)=
\mu_{t_{n}}^{(m)}
\quad\mbox{\rm to the measure}\quad
\Phi_{t_n,t_n}^{(m)}\left(\widetilde{\Phi}^{(m)}_{t_0,t_n}(\widetilde{\mu}_{t_{0}}^{(m)}))\right)=\widetilde{\mu}_{t_{n}}^{(m)}
$$
Recalling that $\widetilde{\mu}_{t_0}^{(m)}=\mu_{t_0}^{(m)}$, we find that
\begin{eqnarray*}
\widetilde{\mu}_{t_{n}}^{(m)}-\mu_{t_{n}}^{(m)}
&
=&\sum_{q=1}^n~\left[\Phi^{(m)}_{t_q,t_n}\left(\widetilde{\Phi}^{(m)}_{t_0,t_q}(\widetilde{\mu}_{t_{0}}^{(m)}))\right)
-\Phi^{(m)}_{t_{q-1},t_n}\left(\widetilde{\Phi}^{(m)}_{t_0,t_{q-1}}(\widetilde{\mu}_{t_{0}}^{(m)}))\right)\right]\\
&=&\sum_{q=1}^n~\left[\Phi^{(m)}_{t_q,t_n}\left(
\Phi^{(m)}_{t_{q-1},t_q}\left(\widetilde{\mu}_{t_{q-1}}^{(m)}\right)+\frac{1}{m^2}~\Ra_{t_{q-1},t_{q}}^{(m)}\left(\widetilde{\mu}_{t_{q-1}}^{(m)}\right)
\right)
-\Phi^{(m)}_{t_{q},t_n}\left(\Phi^{(m)}_{t_{q-1},t_q}\left(\widetilde{\mu}_{t_{q-1}}^{(m)}\right)\right)\right]
\end{eqnarray*}

By (\ref{dec-ref1}),
we prove that
\begin{eqnarray*}
\left[\Phi^{(m)}_{t_p,t_n}(\mu)-\Phi^{(m)}_{t_p,t_n}(\nu)\right](f)&=&
\left[\Psi_{G_{t_{p},t_{n}}^{(m)}}(\mu)-\Psi_{G_{t_{p},t_{n}}^{(m)}}(\nu)\right]P_{t_{p},t_{n}}^{(m)}(f)\\
&=&
\frac{\nu\left(G_{t_{p},t_{n}}^{(m)}\right)}{\mu\left(G_{t_{p},t_{n}}^{(m)}\right)}~
\left[\mu-\nu\right]\left(\frac{G_{t_{p},t_{n}}^{(m)}}{\nu\left(G_{t_{p},t_{n}}^{(m)}\right)}\left[P_{t_{p},t_{n}}^{(m)}(f)-
\Psi_{G_{t_{p},t_{n}}^{(m)}}(\nu)P_{t_{p},t_{n}}^{(m)}(f)
\right]\right)
\end{eqnarray*}
This implies that
$$
\left\Vert 
\Phi^{(m)}_{t_p,t_n}(\mu)-\Phi^{(m)}_{t_p,t_n}(\nu)
\right\Vert_{\rm tv}\leq 2~ g_{t_p,t_n}^{(m)}~\beta\left(P_{t_{p},t_{n}}^{(m)}\right)
\left\Vert \mu-\nu\right\Vert_{\rm tv}$$
with the Dobrushin contraction coefficient $\beta\left(P_{t_{p},t_{n}}^{(m)}\right)(\leq 1)$, and the parameters
$$
g_{t_p,t_n}^{(m)}:=\sup_{x,y}{\frac{G_{t_{p},t_{n}}^{(m)}(x)}{G_{t_{p},t_{n}}^{(m)}(y)}}\leq \exp{(2t_n\sup_{t\in[0,t_n]}\|\Ua_t\|)}~.
$$
This yields the rather crude estimates
\begin{eqnarray*}
m~\left\Vert \widetilde{\mu}_{t_{n}}^{(m)}-\mu_{t_{n}}^{(m)}\right\Vert_{\rm tv}&\leq&
2 m^{-1}~\sum_{p=1}^n~g_{t_p,t_n}^{(m)}~\beta\left(P_{t_{p},t_{n}}^{(m)}\right)~\|\Ua_{t_{q-1}}\|^2\\
&\leq& 2 t_n~\exp{(2t_n\sup_{t\in[0,t_n]}\|\Ua_t\|)}~\sup_{t\in[0,t_n]}{\|\Ua_t\|^2}
\end{eqnarray*}
This ends the proof of the theorem.
\cqfd

Working a little harder, under some regularity conditions, the estimates developed in the proof of the theorem
can be used to obtain uniform estimates w.r.t. the time parameter. For instance, in case \textbf{(D)}, under the stability 
conditions (\ref{stab-condition}), the constant $c_{t_n}$ in theorem~\ref{theo-interm} can be chosen so that $\sup_{n}c_{t_n}<\infty$.
We can extend these uniform results to continuous time models, using the stability analysis of continuous Feynman-Kac semigroups developed in~\cite{dm-Toulouse}.

\section{First order decompositions}

\label{first-order-sect}

The main objective of this section is to prove theorem~\ref{theo-bv1-intro}.
In the further development of this section we let $c$, $c_{n}$,
$c_{t_{n}}$, and $c_{t_{n}}(f)$ be, respectively, some universal
constant, and some finite constants that depend on the parameters
$n$, $t_{n}$, and the pair $(t_{n},f)$, with values that may vary
from line to line but do not depend on the parameters $m$ and $N$.
We also assume that $m$ is chosen so that $\|\Va_{t_{n}}\|\leq c_{t_{n}}~m$,
for any $n\geq0$, and $N\geq m$.

\subsection{Continuous time models}

We start with the continuous time case \textbf{(C)} presented on page~\pageref{casC-ref}.
By lemma~\ref{prop-dec} we find that 
\begin{equation}
\EE\left(W_{t_{n},t_{n+1}}^{N}(f)^{2}~|~\xi_{t_{n}}\right)=\mu_{t_{n}}^{N}\left[\Gamma_{L_{t_{n},\mu_{t_{n}}^{N}}^{(m)}}(f,f)\right]-\mu_{t_{n}}^{N}\left(\left(L_{t_{n},\mu_{t_{n}}^{N}}^{(m)}(f)\right)^{2}\right).\label{lrf}
\end{equation}
Using (\ref{ref-Ka}), for any function $f\in D(L)$, such that $L_{t}(f),L_{t}(f^{2})\in\Ca^{1}([t_{n},t_{n+1}],D(L))$,
we find that 
\[
\EE\left(W_{t_{n},t_{n+1}}^{N}(f)^{2}~|~\xi_{t_{n}}\right)=\mu_{t_{n}}^{N}\left[\Gamma_{L_{t_{n},\mu_{t_{n}}^{N}}}(f,f)\right]~\frac{1}{m}+R_{t_{n}}^{N}\left(f\right)~\frac{1}{m^{2}}
\]
with some remainder term $R_{t_{n}}^{N}(f)$ such that 
\begin{eqnarray*}
\sup_{N\geq1}{\left|R_{t_{n}}^{N}\left(f\right)\right|} & \leq & c_{t_{n}}~\|f^{2}\|_{t_{n}}
\end{eqnarray*}
with the norm $\|f\|_{t_{n}}$ of a function $f\in D(L)$, such that
$L_{t}(f)\in\Ca^{1}([t_{n},t_{n+1}],D(L))$ defined in (\ref{normftn})
By (\ref{Q-ref}), we have 
\[
Q_{t_{p},t_{n}}^{(m)}(f)(x)=\EE\left(f(\Xa_{t_{n}})~\exp{\left[\int_{t_{p}}^{t_{n}}\Va_{\underline{\tau}(s)}(\Xa_{\underline{\tau}(s)})ds\right]}~\left|~\Xa_{t_{p}}=x\right.\right)
\]
with $\underline{\tau}(s)=\sum_{n\geq0}1_{[t_{n},t_{n+1}[}(s)~t_{n}$.
Thus, in case \textbf{(C)}, for any function $f\in D(L)$, the mappings
\[
t\in[t_{p-1},t_{p}]\mapsto L_{t}\left(Q_{t_{p},t_{n}}^{(m)}(f)^{2}\right)\quad\mbox{{\rm and}}\quad t\in[t_{p-1},t_{p}]\mapsto L_{t}^{2}\left(Q_{t_{p},t_{n}}^{(m)}(f)^{2}\right)
\]
are uniformly bounded w.r.t. the parameter $m$, and differentiable
with uniformly bounded derivatives w.r.t. the parameter $m$.

The first assertion of theorem~\ref{theo-bv1-intro} is based on
the first order decompositions of the fluctuation of $\mu_{t_{n}}^{N}$
around its limiting value $\mu_{t_{n}}^{(m)}$ developed in~\cite{dhw2012}.
Using theorem 6.2 in~\cite{dhw2012}, we prove the following proposition.
\begin{prop}\label{prop-refC} For any $N\geq1$ and any $n\in\NN$,
we have 
\begin{eqnarray*}
\sqrt{N}\left[\mu_{t_{n}}^{N}-\mu_{t_{n}}^{(m)}\right] & = & \sum_{p=0}^{n}~\frac{1}{\mu_{t_{p}}^{N}\left(\overline{G}_{t_{p},t_{n}}^{(m,N)}\right)}~W_{t_{p-1},t_{p}}^{N}\left(D_{t_{p},t_{n}}^{(m,N)}(f)\right)\\
 & = & \sum_{p=0}^{n}~~W_{t_{p-1},t_{p}}^{N}\left(D_{t_{p},t_{n}}^{(m,N)}(f)\right)+\frac{1}{\sqrt{N}}~\Ra_{t_{n}}^{(m,N)}(f)
\end{eqnarray*}
with the first order integral operator 
\[
D_{t_{p},t_{n}}^{(m,N)}(f):=\overline{G}_{t_{p},t_{n}}^{(m,N)}\left(P_{t_{p},t_{n}}^{(m)}(f)-\Phi_{t_{p},t_{n}}(\mu_{t_{p-1}}^{N})(f)\right)\quad\mbox{{with}}\quad\overline{G}_{t_{p},t_{n}}^{(m,N)}=G_{t_{p},t_{n}}^{(m)}/\Phi_{t_{p}}^{(m)}(\mu_{t_{p-1}}^{N})\left(G_{t_{p},t_{n}}^{(m)}\right)
\]
and the remainder second order term 
\[
\Ra_{t_{n}}^{(m,N)}(f)=-\sum_{p=0}^{n}~\frac{1}{\mu_{t_{p}}^{N}\left(\overline{G}_{t_{p},t_{n}}^{(m,N)}\right)}W_{t_{p-1},t_{p}}^{N}\left(\overline{G}_{t_{p},t_{n}}^{(m,N)}\right)~W_{t_{p-1},t_{p}}^{N}\left(D_{t_{p},t_{n}}^{(m,N)}(f)\right)
\]
In the above, we have used the convention $W_{t_{-1},t_{0}}^{N}=\sqrt{N}[\mu_{t_{0}}^{N}-\mu_{t_{0}}]$,
for $p=0$. \end{prop}

We are now in position to prove the bias and the variance estimates
stated in theorem~\ref{theo-bv1-intro} for the continuous time models.

\textbf{Proof of theorem~\ref{theo-bv1-intro} - case (C):}

Firstly, the first order decomposition stated above clearly implies
that 
\[
N~\EE\left(\left[\mu_{t_{n}}^{N}-\mu_{t_{n}}^{(m)}\right]\right)=\EE\left(\Ra_{t_{n}}^{(m,N)}(f)\right).
\]
On the other hand, we have 
\[
\left|\EE\left(\Ra_{t_{n}}^{(m,N)}(f)\right)\right|\leq c_{t_{n}}^{1}~\sum_{p=0}^{n}~\EE\left(W_{t_{p-1},t_{p}}^{N}\left(\overline{G}_{t_{p},t_{n}}^{(m,N)}\right)^{2}\right)^{1/2}\EE\left(W_{t_{p-1},t_{p}}^{N}\left(D_{t_{p},t_{n}}^{(m,N)}(f)\right)^{2}\right)^{1/2}.
\]
To get one step further, we use the fact that 
\[
\begin{array}{l}
\EE\left(W_{t_{p-1},t_{p}}^{N}\left(D_{t_{p},t_{n}}^{(m,N)}(f)\right)^{2}\right)\\
\\
=\EE\left(\mu_{t_{p-1}}^{N}\left[\Gamma_{L_{t_{p-1},\mu_{t_{p-1}}^{N}}}(D_{t_{p},t_{n}}^{(m,N)}(f),D_{t_{p},t_{n}}^{(m,N)}(f))\right]\right)~{\displaystyle \frac{1}{m}+\EE\left(R_{t_{p-1}}^{N}\left(D_{t_{p},t_{n}}^{(m,N)}(f)\right)\right)~{\displaystyle \frac{1}{m^{2}}}}.
\end{array}
\]

After some elementary manipulations we prove that 
\begin{equation}
\sup_{0\leq p\leq n}{\EE\left(W_{t_{p-1},t_{p}}^{N}\left(D_{t_{p},t_{n}}^{(m,N)}(f)\right)^{2}\right)}\leq c_{t_{n}}(f)/{m}\quad\mbox{{\rm and}}\quad\sup_{0\leq p\leq n}{\EE\left(W_{t_{p-1},t_{p}}^{N}\left(\overline{G}_{t_{p},t_{n}}^{(m,N)}\right)^{2}\right)}\leq c_{t_{n}}/{m}.\label{refW4}
\end{equation}
This ends the proof of the bias estimate.

The proof of the variance estimates is based on the following technical
lemma. \begin{lem}\label{tech-lem} For any $f$ with $\mbox{{\rm osc}}(f)\leq1$,
we have the fourth conditional moment estimate 
\begin{eqnarray*}
\EE\left(W_{t_{n},t_{n+1}}^{N}(f)^{4}~|~\xi_{t_{n}}\right) & \leq & \frac{1}{N}~\EE\left(W_{t_{n},t_{n+1}}^{N}(f)^{2}~|~\xi_{t_{n}}\right)+6~\EE\left(W_{t_{n},t_{n+1}}^{N}(f)^{2}~|~\xi_{t_{n}}\right)^{2}.
\end{eqnarray*}
\end{lem} \proof With some elementary computation, we have that
\[
\begin{array}{l}
\EE\left(W_{t_{n},t_{n+1}}^{N}(f)^{4}~|~\xi_{t_{n}}\right)\\
\\
=\frac{1}{N}\int\mu_{t_{n}}^{N}(du)~\Ka_{t_{n},t_{n+1},\mu_{t_{n}}^{N}}(u,dv)~\left(f(v)-\Ka_{t_{n},t_{n+1},\mu_{t_{n}}^{N}}(f)(u)\right)^{4}\\
\\
+6\left(1-\frac{1}{N}\right)~\int_{u\not=u^{\prime}}\mu_{t_{n}}^{N}(du)\mu_{t_{n}}^{N}(du^{\prime})~\Ka_{t_{n},t_{n+1},\mu_{t_{n}}^{N}}(u,dv)~\left(f(v)-\Ka_{t_{n},t_{n+1},\mu_{t_{n}}^{N}}(f)(u)\right)^{2}\\
\\
\hskip7cm\Ka_{t_{n},t_{n+1},\mu_{t_{n}}^{N}}(u^{\prime},dv^{\prime})~\left(f(v^{\prime})-\Ka_{t_{n},t_{n+1},\mu_{t_{n}}^{N}}(f)(u^{\prime})\right)^{2}.
\end{array}
\]
This implies that 
\begin{eqnarray*}
\EE\left(W_{t_{n},t_{n+1}}^{N}(f)^{4}~|~\xi_{t_{n}}\right) & \leq & \frac{1}{N}\int\mu_{t_{n}}^{N}(du)~\Ka_{t_{n},t_{n+1},\mu_{t_{n}}^{N}}(u,dv)~\left(f(v)-\Ka_{t_{n},t_{n+1},\mu_{t_{n}}^{N}}(f)(u)\right)^{4}\\
 &  & +6\left[\int\mu_{t_{n}}^{N}(du)\Ka_{t_{n},t_{n+1},\mu_{t_{n}}^{N}}(u,dv)~\left(f(v)-\Ka_{t_{n},t_{n+1},\mu_{t_{n}}^{N}}(f)(u)\right)^{2}\right]^{2}.
\end{eqnarray*}
The end of the proof is based on the fact that 
\[
\left(f(v)-\Ka_{t_{n},t_{n+1},\mu_{t_{n}}^{N}}(f)(u)\right)^{4}\leq\left(f(v)-\Ka_{t_{n},t_{n+1},\mu_{t_{n}}^{N}}(f)(u)\right)^{2}
\]

as soon as $\mbox{{\rm osc}}(f)\leq1$. This ends the proof of the
lemma. \cqfd

Combining this lemma with (\ref{refW4}), we prove that 
\[
\sup_{0\leq p\leq n}{\EE\left(W_{t_{p-1},t_{p}}^{N}\left(D_{t_{p},t_{n}}^{(m,N)}(f)\right)^{4}\right)}\leq c_{t_{n}}(f)~\frac{1}{m}~\left[\frac{1}{N}+\frac{1}{m}\right]
\]
and 
\[
\sup_{0\leq p\leq n}{\EE\left(W_{t_{p-1},t_{p}}^{N}\left(\overline{G}_{t_{p},t_{n}}^{(m,N)}\right)^{4}\right)}\leq c_{t_{n}}(f)~\frac{1}{m}~\left[\frac{1}{N}+\frac{1}{m}\right].
\]
This implies that 
\[
\EE\left(\Ra_{t_{n}}^{(m,N)}(f)^{2}\right)^{1/2}\leq c_{t_{n}}~\sum_{p=0}^{n}~\EE\left[W_{t_{p-1},t_{p}}^{N}\left(\overline{G}_{t_{p},t_{n}}^{(m,N)}\right)^{4}\right]^{1/4}~\EE\left[W_{t_{p-1},t_{p}}^{N}\left(D_{t_{p},t_{n}}^{(m,N)}(f)\right)^{4}\right]^{1/4}
\]
and therefore 
\[
\EE\left(\Ra_{t_{n}}^{(m,N)}(f)^{2}\right)\leq c_{t_{n}}(f)~(1+m/N).
\]
Using the fact that 
\[
N~\EE\left(\left(\mu_{t_{n}}^{N}(f)-\mu_{t_{n}}^{(m)}(f)\right)^{2}\right)\leq2~\left[\EE\left(\left(\sum_{p=0}^{n}~~W_{t_{p-1},t_{p}}^{N}\left(D_{t_{p},t_{n}}^{(m,N)}(f)\right)\right)^{2}\right)+\frac{1}{N}~\EE\left(\Ra_{t_{n}}^{(m,N)}(f)^{2}\right)\right]
\]
with 
\begin{eqnarray*}
\EE\left(\left(\sum_{p=0}^{n}~~W_{t_{p-1},t_{p}}^{N}\left(D_{t_{p},t_{n}}^{(m,N)}(f)\right)\right)^{2}\right) & = & \sum_{p=0}^{n}~\EE\left(\left(W_{t_{p-1},t_{p}}^{N}\left(D_{t_{p},t_{n}}^{(m,N)}(f)\right)\right)^{2}\right)\leq c_{t_{n}}(f)
\end{eqnarray*}
we conclude that 
\[
N~\EE\left(\left(\mu_{t_{n}}^{N}(f)-\mu_{t_{n}}^{(m)}(f)\right)^{2}\right)\leq c_{t_{n}}(f)\left(1+\frac{1}{N}+\frac{m}{N^{2}}\right).
\]
This ends the proof of theorem~\ref{theo-bv1-intro}. \cqfd

\subsection{Discrete time models}

The main objective of this section is to prove theorem~\ref{theo-bv1-intro}
for the discrete time models related to case \textbf{(D)}. In this
situation, we recall that $\mu_{t_{nm}}=\eta_{n}$, for any integer
parameter $n\in\NN$. As usual, the approximation measures $\mu_{t_{n}}^{N}$
are defined as the occupation measures of the mean field particle
model of the McKean distribution flow $\mu_{t_{n}}$. By lemma~\ref{lem-ref-D},
for any $k\in\NN$ we have 
\begin{eqnarray*}
L_{k-\frac{1}{m},\mu_{k-\frac{1}{m}}}^{(m)} & := & \Ka_{k-\frac{1}{m},~k,~\mu_{k-\frac{1}{m}}}-Id=\Sa_{k-\frac{1}{m},\mu_{k-\frac{1}{m}}}M_{k}-Id
\end{eqnarray*}
and for any $(k-1)m<p<km$, with $k\in\NN$, we have 
\[
\Ma_{t_{p-1},t_{p}}=Id\Longrightarrow L_{t_{p-1},\mu_{t_{p-1}}}^{(m)}=\Sa_{t_{p-1},\mu_{t_{p-1}}}-Id.
\]
When the Markov transport equation (\ref{MT}) is met for some Markov
transitions $\Sa_{t_{n},\mu}$ satisfying the first order decomposition
(\ref{1st-exp-ref}), we have for any $k\in\NN$ and any $(k-1)m<p<km$
\[
L_{t_{p-1},\mu_{t_{p-1}}}^{(m)}=\widehat{L}_{t_{p-1},\mu_{t_{p-1}}}~\frac{1}{m}+\widehat{R}_{t_{p-1},\mu_{t_{p-1}}}~\frac{1}{m^{2}}.
\]
Using the same line of argument as in (\ref{lrf}), we prove that
\begin{eqnarray*}
\EE\left(W_{t_{p-1},t_{p}}^{N}(f)^{2}~|~\xi_{t_{n}}\right) & = & \mu_{t_{p-1}}^{N}\left[\Gamma_{L_{t_{p-1},\mu_{t_{p-1}}^{N}}^{(m)}}(f,f)\right]-\mu_{t_{p-1}}^{N}\left(\left(L_{t_{p-1},\mu_{t_{p-1}}^{N}}^{(m)}(f)\right)^{2}\right)\\
 & = & \mu_{t_{p-1}}^{N}\left[\Gamma_{\widehat{L}_{t_{p-1},\mu_{t_{p-1}}^{N}}}(f,f)\right]~\frac{1}{m}+R_{t_{p-1}}^{N}\left(f\right)~\frac{1}{m^{2}}
\end{eqnarray*}
with some remainder term $R_{t_{p-1}}^{N}(f)$ such that 
\[
\sup_{N\geq1}{\left|R_{t_{p-1}}^{N}\left(f\right)\right|}\leq c~\left\Vert \Ua_{t_{p-1}}\right\Vert ^{2}~\mbox{{\rm osc}}(f)^{2}.
\]
This clearly implies that 
\begin{equation}
\EE\left(W_{t_{p-1},t_{p}}^{N}(f)^{2}~|~\xi_{t_{n}}\right)\leq c~\left\Vert \Ua_{t_{p-1}}\right\Vert ^{2}~\mbox{{\rm osc}}(f)^{2}~\frac{1}{m}\quad\mbox{{\rm and}}\quad\EE\left(W_{k-\frac{1}{m},k}^{N}(f)^{2}~|~\xi_{k-\frac{1}{m}}\right)\leq\mbox{{\rm osc}}(f)^{2}.\label{inter-1}
\end{equation}
As in proposition~\ref{prop-refC}, we prove the following decomposition.
\begin{prop} For any $N\geq1$, $f\in\Ba_{b}(E)$, and any $n\in\NN$
we have the decomposition 
\[
[\mu_{n}^{N}-\eta_{n}](f)=A_{n}^{N}+B_{n}^{N}
\]
with 
\begin{eqnarray*}
A_{n}^{N} & := & \sum_{1\leq k\leq n}\sum_{p=(k-1)m+1}^{(k-1)m+(m-1)}\frac{1}{\mu_{t_{p}}^{N}\left(\overline{G}_{t_{p},n}^{(m,N)}\right)}~W_{t_{p-1},t_{p}}^{N}\left(D_{t_{p},n}^{(m,N)}(f)\right)
\end{eqnarray*}
and 
\[
B_{n}^{N}=\sum_{0\leq k\leq n}\frac{1}{\mu_{k}^{N}\left(\overline{G}_{k,n}^{(m,N)}\right)}~W_{{k-\frac{1}{m}},k}^{N}\left(D_{k,n}^{(m,N)}(f)\right).
\]
In the above, the function $\overline{G}_{t_{p},n}^{(m,N)}$ and the
integral operator $D_{t_{p},n}^{(m,N)}(f)$ are defined in proposition~\ref{prop-refC},
and for $k=0$, we have used the convention $W_{{-\frac{1}{m}},0}^{N}=W_{0}^{N}$.
\end{prop}

To analyze the bias and the variance, we also need to consider the
first order decompositions presented in the following corollary. \begin{cor}
For any $N\geq1$ and any $n\in\NN$ we have the decomposition 
\[
A_{n}^{N}=A_{n}^{(N,1)}+\frac{1}{\sqrt{N}}~A_{n}^{(N,2)}\quad\mbox{{\rm and}}\quad B_{n}^{N}=B_{n}^{(N,1)}+\frac{1}{\sqrt{N}}~B_{n}^{(N,2)}
\]
with the first order terms 
\begin{eqnarray*}
A_{n}^{(N,1)} & = & \sum_{1\leq k\leq n}\sum_{p=(k-1)m+1}^{(k-1)m+(m-1)}~W_{t_{p-1},t_{p}}^{N}\left(D_{t_{p},n}^{(m,N)}(f)\right)\\
B_{n}^{(N,1)} & = & \sum_{0\leq k\leq n}W_{{k-\frac{1}{m}},k}^{N}\left(D_{k,n}^{(m,N)}(f)\right)
\end{eqnarray*}
and the remainder second order terms 
\begin{eqnarray*}
A_{n}^{(N,2)} & = & -\sum_{1\leq k\leq n}\sum_{p=(k-1)m+1}^{(k-1)m+(m-1)}\frac{1}{\mu_{t_{p}}^{N}\left(\overline{G}_{t_{p},n}^{(m,N)}\right)}~W_{t_{p-1},t_{p}}^{N}\left(\overline{G}_{t_{p},n}^{(m,N)}\right)~W_{t_{p-1},t_{p}}^{N}\left(D_{t_{p},n}^{(m,N)}(f)\right)\\
B_{n}^{(N,2)} & = & -\sum_{0\leq k\leq n}\frac{1}{\mu_{k}^{N}\left(\overline{G}_{k,n}^{(m,N)}\right)}~W_{{k-\frac{1}{m}},k}^{N}\left(\overline{G}_{k,n}^{(m,N)}\right)~W_{{k-\frac{1}{m}},k}^{N}\left(D_{k,n}^{(m,N)}(f)\right)
\end{eqnarray*}
\end{cor}

Now we come to the proof of the bias and the variance estimates presented
in theorem~\ref{theo-bv1-intro}.

In the further development of this section $f$ stands for some bounded
function s.t. $\mbox{{\rm osc}}(f)\leq1$.

By construction, the random fields $W_{t_{p-1},t_{p}}^{N}$ and $W_{t_{q-1},t_{q}}^{N}$
are uncorrelated for any $p\not=q$. Combining this property with
the estimates (\ref{inter-1}) we prove that 
\begin{eqnarray*}
\EE\left(\left(A_{n}^{(N,1)}\right)^{2}\right) & = & \sum_{1\leq k\leq n}\sum_{p=(k-1)m+1}^{(k-1)m+(m-1)}~\EE\left(\left[W_{t_{p-1},t_{p}}^{N}\left(D_{t_{p},n}^{(m,N)}(f)\right)\right]^{2}\right).
\end{eqnarray*}
On the other hand, we have 
\[
\EE\left(\left[W_{t_{p-1},t_{p}}^{N}\left(D_{t_{p},n}^{(m,N)}(f)\right)\right]^{2}\right)\leq\frac{c}{m}~\left\Vert \Ua_{t_{p-1}}\right\Vert ^{2}~\left(g_{t_{p},n}^{(m)}~\mbox{{\rm osc}}\left(P_{t_{p},n}^{(m)}(f)\right)\right)^{2}
\]
with 
\[
g_{t_{p},n}^{(m)}:=\sup_{x,y}{\left[{G_{t_{p},n}^{(m)}(x)}/{G_{t_{p},n}^{(m)}(y)}\right]}.
\]
We prove the last assertion using the fact that 
\[
\mbox{{\rm osc}}\left(D_{t_{p},n}^{(m,N)}(f)\right)/2\leq\left\Vert D_{t_{p},n}^{(m,N)}(f)\right\Vert \leq g_{t_{p},n}^{(m)}~\mbox{{\rm osc}}\left(P_{t_{p},n}^{(m)}(f)\right).
\]
By the semigroup formulae (\ref{sg-formula}), for any $p=(k-1)m+r$,
with $r<m$ we find that 
\[
P_{t_{p},n}^{(m)}=P_{(k-1),n}\quad\mbox{{\rm and}}\quad g_{t_{p},n}^{(m)}\leq g_{k-1,n}\times\sup_{x,y}{\left(\frac{G_{k-1}(x)}{G_{k-1}(y)}\right)^{1/m}}=g_{k-1,n}\times g_{k-1,k}^{1/m}
\]
with 
\[
g_{k-1,n}:=\sup_{x,y}{\left[Q_{(k-1),n}(1)(x)/Q_{(k-1),n}(1)(y)\right]}.
\]

This implies that 
\[
\EE\left(\left(A_{n}^{(N,1)}\right)^{2}\right)\leq c~\sum_{1\leq k\leq n}\left\Vert \log{G_{k-1}}\right\Vert ^{2}~\left(g_{k-1,n}\times g_{k-1,k}^{1/m}~\mbox{{\rm osc}}\left(P_{(k-1),n}(f)\right)\right)^{2}.
\]

In the same way, we prove that 
\begin{eqnarray*}
\EE\left(\left(B_{n}^{(N,1)}\right)^{2}\right) & = & \sum_{0\leq k\leq n}\EE\left(\left[W_{{k-\frac{1}{m}},k}^{N}\left(D_{k,n}^{(m,N)}(f)\right)\right]^{2}\right)\\
 & \leq & c~\sum_{0\leq k\leq n}~\left[g_{k,n}^{(m)}~\mbox{{\rm osc}}\left(P_{k,n}^{(m)}(f)\right)\right]^{2}=c~\sum_{0\leq k\leq n}~\left[g_{k,n}~\mbox{{\rm osc}}\left(P_{k,n}(f)\right)\right]^{2}.
\end{eqnarray*}

On the other hand, we have 
\begin{equation}
\begin{array}{l}
\left|\EE\left(A_{n}^{(N,2)}\right)\right|\\
\\
\leq{\displaystyle \sum_{1\leq k\leq n}g_{k-1,n}g_{k-1,k}^{1/m}~\sum_{p=(k-1)m+1}^{(k-1)m+(m-1)}}\\
\\
\hskip4cm\EE\left(\left(W_{t_{p-1},t_{p}}^{N}\left(\overline{G}_{t_{p},n}^{(m,N)}\right)\right)^{2}\right)^{1/2}\EE\left(\left(W_{t_{p-1},t_{p}}^{N}\left(D_{t_{p},n}^{(m,N)}(f)\right)\right)^{2}\right)^{1/2}\\
\\
\leq{\displaystyle \sum_{0\leq k<n}\left\Vert \log{G_{k}}\right\Vert ^{2}~\left(g_{k,n}g_{k,k+1}\right)^{3}~\mbox{{\rm osc}}\left(P_{k,n}(f)\right)}.
\end{array}\label{bias-ref1}
\end{equation}

Arguing as in the end of the proof of theorem~\ref{theo-bv1-intro},
we can also check that 
\[
\begin{array}{l}
\EE\left(\left(A_{n}^{(N,2)}\right)^{2}\right)^{1/2}\\
\\
\leq~{\displaystyle \sum_{1\leq k\leq n}~g_{k-1,n}g_{k-1,k}^{1/m}\sum_{p=(k-1)m+1}^{(k-1)m+(m-1)}\EE\left(\left(W_{t_{p-1},t_{p}}^{N}\left(\overline{G}_{t_{p},n}^{(m,N)}\right)~W_{t_{p-1},t_{p}}^{N}\left(D_{t_{p},n}^{(m,N)}(f)\right)\right)^{2}\right)^{1/2}}\\
\\
\\
\leq{\displaystyle \sum_{1\leq k\leq n}g_{k-1,n}g_{k-1,k}^{1/m}~\sum_{p=(k-1)m+1}^{(k-1)m+(m-1)}}\\
\\
\hskip4cm\EE\left(\left(W_{t_{p-1},t_{p}}^{N}\left(\overline{G}_{t_{p},n}^{(m,N)}\right)\right)^{4}\right)^{1/4}\EE\left(\left(W_{t_{p-1},t_{p}}^{N}\left(D_{t_{p},n}^{(m,N)}(f)\right)\right)^{4}\right)^{1/4}.
\end{array}
\]

On the other hand, using lemma~\ref{tech-lem} we have 
\[
\begin{array}{l}
\EE\left(\left[W_{t_{p-1},t_{p}}^{N}\left(\overline{G}_{t_{p},n}^{(m,N)}\right)\right]^{4}\right)\\
\\
\leq\frac{1}{N}~(2g_{t_{p},n}^{(m)})^{2}~\EE\left(\left[W_{t_{p-1},t_{p}}^{N}\left(\overline{G}_{t_{p},n}^{(m,N)}\right)\right]^{2}\right)+6~\EE\left(\left[W_{t_{p-1},t_{p}}^{N}\left(\overline{G}_{t_{p},n}^{(m,N)}\right)\right]^{2}\right)^{2}\\
\\
\leq c~(g_{t_{p},n}^{(m)})^{4}~\left(\frac{1}{N}~\|\Ua_{t_{p-1}}\|^{2}~\frac{1}{m}+\|\Ua_{t_{p-1}}\|^{4}~\frac{1}{m^{2}}\right)
\end{array}
\]
from which we find the crude upper bound 
\[
\EE\left(\left[W_{t_{p-1},t_{p}}^{N}\left(\overline{G}_{t_{p},n}^{(m,N)}\right)\right]^{4}\right)\leq c~g_{k-1,n}^{4}g_{k-1,k}^{4}~\left(\|\log{G_{k-1}}\|\vee1\right)^{4} /{m^{2}}.
\]
for any $N\geq m$, and for any $p=(k-1)m+r$, with $r<m$. Using
the same line of arguments, we have 
\[
\begin{array}{l}
\EE\left(\left[W_{t_{p-1},t_{p}}^{N}\left(D_{t_{p},n}^{(m,N)}(f)\right)\right]^{4}\right)\\
\\
\leq\frac{1}{N}~~\left(2~g_{t_{p},n}^{(m)}~\mbox{{\rm osc}}\left(P_{t_{p},n}^{(m)}(f)\right)\right)^{2}\EE\left(\left[W_{t_{p-1},t_{p}}^{N}\left(D_{t_{p},n}^{(m,N)}(f)\right)\right]^{2}\right)+6~\EE\left(\left[W_{t_{p-1},t_{p}}^{N}\left(D_{t_{p},n}^{(m,N)}(f)\right)\right]^{2}\right)^{2}\\
\\
\leq c~g_{k-1,n}^{4}g_{k-1,k}^{4/m}~{\displaystyle \mbox{{\rm osc}}\left(P_{(k-1),n}(f)\right)^{4}~\left(\|\log{G_{k-1}}\|\vee1\right)^{4} /{m^{2}}}
\end{array}
\]
for any $N\geq m$. This implies that 
\[
\EE\left(\left(A_{n}^{(N,2)}\right)^{2}\right)^{1/2}\leq c~{\displaystyle \sum_{1\leq k\leq n}g_{k-1,n}^{3}g_{k-1,k}^{3}\left(\|\log{G_{k-1}}\|\vee1\right)^{2}{\displaystyle \mbox{{\rm osc}}\left(P_{(k-1),n}(f)\right)}}
\]
for any $N\geq m$. One concludes that 
\[
\EE\left(\left(A_{n}^{N}\right)^{2}\right)\leq c~\left(a_{1,n}+\frac{1}{N}~a_{2,n}^{2}\right)\leq c~a_{2,n}~\left(1+\frac{1}{N}~a_{2,n}\right)
\]
with 
\[
a_{1,n}:=~\sum_{0\leq k<n}\left\Vert \log{G_{k}}\right\Vert ^{2}~\left(g_{k,n}g_{k,k+1}~\mbox{{\rm osc}}\left(P_{k,n}(f)\right)\right)^{2}\leq a_{2,n}
\]
and 
\[
a_{2,n}:={\displaystyle \sum_{0\leq k<n}g_{k,n}^{3}g_{k,k+1}^{3}\left(\|\log{G_{k}}\|\vee1\right)^{2}{\displaystyle \mbox{{\rm osc}}\left(P_{k,n}(f)\right)}}.
\]
In much the same way, we prove that 
\[
\EE\left(\left(B_{n}^{(N,2)}\right)^{2}\right)^{1/2}\leq\sum_{0\leq k\leq n}g_{k,n}~\EE\left(\left[W_{{k-\frac{1}{m}},k}^{N}\left(\overline{G}_{k,n}^{(m,N)}\right)\right]^{4}\right)^{1/4}\EE\left(\left[W_{{k-\frac{1}{m}},k}^{N}\left(D_{k,n}^{(m,N)}(f)\right)\right]^{4}\right)^{1/4}.
\]
Using the fact that 
\[
\EE\left(\left[W_{{k-\frac{1}{m}},k}^{N}\left(\overline{G}_{k,n}^{(m,N)}\right)\right]^{2}\right)\leq c~g_{k,n}^{2}
\]
and 
\[
\EE\left(\left[W_{{k-\frac{1}{m}},k}^{N}\left(D_{k,n}^{(m,N)}(f)\right)\right]^{2}\right)\leq c~g_{k,n}^{2}~\mbox{{\rm osc}}\left(P_{k,n}(f)\right)^{2}
\]
we check that 
\begin{eqnarray}
\left|\EE\left(B_{n}^{(N,2)}\right)\right| & \leq & \sum_{0\leq k\leq n}g_{k,n}~\EE\left(\left[W_{{k-\frac{1}{m}},k}^{N}\left(\overline{G}_{k,n}^{(m,N)}\right)\right]^{2}\right)^{1/2}\EE\left(\left[W_{{k-\frac{1}{m}},k}^{N}\left(D_{k,n}^{(m,N)}(f)\right)\right]^{2}\right)^{1/2}\nonumber \\
 & \leq & c~\sum_{0\leq k\leq n}~g_{k,n}^{3}~\mbox{{\rm osc}}\left(P_{k,n}(f)\right).\label{bias-ref2}
\end{eqnarray}
Combining (\ref{bias-ref1}) with (\ref{bias-ref2}), we obtain the
bias estimate 
\[
\begin{array}{l}
N~\left|\EE\left(\mu_{n}^{N}(f)\right)-\eta_{n}(f)\right|\\
\\
\leq c~\left(\sum_{0\leq k\leq n}~g_{k,n}^{3}~\mbox{{\rm osc}}\left(P_{k,n}(f)\right)+{\displaystyle \sum_{0\leq k<n}\left\Vert \log{G_{k}}\right\Vert ^{2}~\left(g_{k,n}g_{k,k+1}\right)^{3}~\mbox{{\rm osc}}\left(P_{k,n}(f)\right)}\right)\\
\\
\leq c~{\displaystyle \sum_{0\leq k\leq n}\left(\left\Vert \log{G_{k}}\right\Vert \vee1\right)^{2}~\left(g_{k,n}g_{k,k+1}\right)^{3}~\mbox{{\rm osc}}\left(P_{k,n}(f)\right)}.
\end{array}
\]

In addition, we have that 
\[
\begin{array}{l}
\EE\left(\left[W_{{k-\frac{1}{m}},k}^{N}\left(\overline{G}_{k,n}^{(m,N)}\right)\right]^{4}\right)\\
\\
\leq c~\left[{\displaystyle \frac{g_{k,n}^{2}}{N}~\EE\left(\left[W_{{k-\frac{1}{m}},k}^{N}\left(\overline{G}_{k,n}^{(m,N)}\right)\right]^{2}\right)+\EE\left(\left[W_{{k-\frac{1}{m}},k}^{N}\left(\overline{G}_{k,n}^{(m,N)}\right)\right]^{2}\right)^{2}}\right]\leq c~g_{k,n}^{4}
\end{array}
\]
and 
\[
\begin{array}{l}
\EE\left(\left[W_{{k-\frac{1}{m}},k}^{N}\left(D_{k,n}^{(m,N)}(f)\right)\right]^{4}\right)\\
\\
\leq c~\left[\displaystyle\frac{g_{k,n}^{2}}{N}~\mbox{{\rm osc}}\left(P_{k,n}(f)\right)^{2}~\EE\left(\left[W_{{k-\frac{1}{m}},k}^{N}\left(D_{k,n}^{(m,N)}(f)\right)\right]^{2}\right)+\EE\left(\left[W_{{k-\frac{1}{m}},k}^{N}\left(D_{k,n}^{(m,N)}(f)\right)\right]^{2}\right)^{2}\right]\\
\\
\leq c~g_{k,n}^{4}\mbox{{\rm osc}}\left(P_{k,n}(f)\right)^{4}.
\end{array}
\]
This implies that 
\[
\EE\left(\left(B_{n}^{(N,2)}\right)^{2}\right)^{1/2}\leq c~\sum_{0\leq k\leq n}g_{k,n}^{3}~\mbox{{\rm osc}}\left(P_{k,n}(f)\right).
\]
We conclude that 
\[
\EE\left(\left(B_{n}^{N}\right)^{2}\right)\leq c~\left(b_{1,n}+\frac{1}{N}~b_{2,n}^{2}\right)\leq c~b_{2,n}~\left(1+\frac{1}{N}~b_{2,n}\right)
\]
with 
\[
b_{1,n}:=\sum_{0\leq k\leq n}~\left[g_{k,n}~\mbox{{\rm osc}}\left(P_{k,n}(f)\right)\right]^{2}\leq b_{2,n}\leq a_{2,n}\quad\mbox{{\rm and}}\quad b_{2,n}:=\sum_{0\leq k\leq n}g_{k,n}^{3}~\mbox{{\rm osc}}\left(P_{k,n}(f)\right).
\]
This yields the variance estimate 
\[
N~\EE\left(\left(\mu_{n}^{N}(f)-\eta_{n}(f)\right)^{2}\right)\leq c~a_{2,n}~\left(1+\frac{1}{N}~a_{2,n}\right).
\]
This ends the proof of the theorem.\cqfd

\section{Appendix}

\subsection{Proof of theorem~\ref{theo1}}

We let $\underline{\tau}$ and $\overline{\tau}$ be the mappings
on $\RR_{+}$ defined by 
\[
\underline{\tau}(s)=\sum_{n\geq0}1_{[t_{n},t_{n+1}[}(s)~t_{n}\quad\mbox{{\rm and}}\quad\overline{\tau}(s)=\sum_{n\geq0}1_{[t_{n},t_{n+1}[}(s)~t_{n+1}.
\]
With this notation, we clearly have that 
\[
\int_{0}^{t_{n}}\Va_{\underline{\tau}(s)}(\Xa_{\underline{\tau}(s)})ds=\sum_{0\leq p<n}\int_{t_{p}}^{t_{p+1}}\Va_{\underline{\tau}(s)}(\Xa_{\underline{\tau}(s)})ds=\sum_{0\leq p<n}\Va_{t_{p}}(\Xa_{t_{p}})/m.
\]

In case \textbf{(D)} we readily check that $\gamma_{\lfloor t\rfloor}=\nu_{\lfloor t\rfloor}$,
for any $t\in\RR_{+}$. More precisely, we have 
\[
n=km+r\quad\mbox{{\rm with}}\quad k\geq0,~~0\leq r<m\Longrightarrow t_{n}=k+r/m
\]
and the pair of processes $(\Va_{s},\Xa_{s})$ only change at integer
times, that is we have that 
\[
t_{pm}=p\leq s<t_{(p+1)m}=(p+1)\Longrightarrow\Va_{s}=\log{G_{p}}\quad\mbox{{\rm and}}\quad\Xa_{s}=X_{p}
\]
so that for any $q\in\NN$ we have that 
\[
t_{q}\leq s<t_{q+1}\Longrightarrow\Va_{s}=\Va_{t_{q}}\quad\mbox{{\rm and}}\quad\Xa_{s}=\Xa_{t_{q}}.
\]
Furthermore, using the fact that 
\[
\int_{0}^{t_{n}}\Va_{s}(\Xa_{s})ds=\sum_{0\leq p<n}\Va_{t_{p}}(\Xa_{t_{p}})/m
\]
we readily check that 
\[
\nu_{t_{n}}(f)=\EE\left(f(\Xa_{t_{n}})~\prod_{0\leq p<n}e^{\Va_{t_{p}}(\Xa_{t_{p}})/m}\right)\qquad\nu_{t_{nm}}=\gamma_{n}\quad\mbox{{\rm and}}\quad\mu_{t_{nm}}=\eta_{n}.
\]

Now we come to the case \textbf{(C)}. We have the Ito formula 
\[
d\Va_{t}(\Xa_{t})=\left(\frac{\partial}{\partial t}+L_{t}\right)(\Va_{t})(\Xa_{t})+dM_{t}(\Va)
\]
with a martingale term $M_{t}(\Va)$ with predictable angle bracket
\[
d\langle M(\Va)\rangle_{t}=\Gamma_{L_{t}}(\Va_{t},\Va_{t})(\Xa_{t})dt
\]
defined in terms of the carré du champ $\Gamma_{L_{t}}$ operator
associated with the generator $L_{t}$ and defined for any $f\in D(L)$
by the following formula 
\[
\Gamma_{L_{t}}(f,f)(x)=L_{t}((f-f(x))^{2})(x)=L_{t}(f^{2})(x)-2f(x)L_{t}(f)(x).
\]
We recall that the predictable process $\langle M(\Va)\rangle_{t}$
is the unique right-continuous and increasing predictable such that
the random process $M_{t}(\Va)^{2}-\langle M(\Va)\rangle_{t}$ is
again a martingale. We also recall that $M_{t}(\Va)^{2}-\left[M(\Va)\right]_{t}$
is also a martingale, for the quadratic variation $\left[M(\Va)\right]$
of the process $M(\Va)$ defined as 
\[
\left[M(\Va)\right]_{t}=\lim_{\|\pi\|\rightarrow0}\sum_{k=1}^{n}\left(M_{t_{k}}(\Va)-M_{t_{k-1}}(\Va)\right)^{2}
\]
where $\pi$ ranges over partitions of the interval $[0,t]$, and
the norm $\|\pi\|$ of the partition $\pi$ is the size of the mesh.
For continuous martingales $M_{t}(\Va)$, it is well known that $\left[M(\Va)\right]=\langle M(\Va)\rangle$.

Using an elementary integration by part formula, for any $p\geq0$
we have 
\[
\int_{t_{p}}^{t_{p+1}}\left(\Va_{s}(\Xa_{s})-\Va_{t_{p}}(\Xa_{t_{p}})\right)ds=\int_{t_{p}}^{t_{p+1}}(t_{p+1}-s)~d\Va_{s}(\Xa_{s})
\]
from which we conclude that 
\[
\int_{0}^{t_{n}}\Va_{s}(\Xa_{s})ds=\int_{0}^{t_{n}}\Va_{\underline{\tau}(s)}(\Xa_{\underline{\tau}(s)})ds+\int_{0}^{t_{n}}(\overline{\tau}(s)-s)~d\Va_{s}(\Xa_{s}).
\]
We set 
\[
R_{t_{n}}(f)=\EE\left(f(\Xa_{[0,t_{n}]})~\left[e^{\int_{0}^{t_{n}}\Va_{s}(\Xa_{s})ds}-e^{\int_{0}^{t_{n}}\Va_{\underline{\tau}(s)}(\Xa_{\underline{\tau}(s)})ds}\right]\right).
\]
Using the fact that $|e^{x}-e^{y}|\leq|x-y|~|e^{x}+e^{y}|$, we prove
that 
\[
\left|R_{t_{n}}(f)\right|\leq\EE\left(\left[\int_{0}^{t_{n}}(\overline{\tau}(s)-s)~d\Va_{s}(\Xa_{s})~\right]\left[e^{\int_{0}^{t_{n}}\Va_{s}(\Xa_{s})ds}+e^{\int_{0}^{t_{n}}\Va_{\underline{\tau}(s)}(\Xa_{\underline{\tau}(s)})ds}\right]\right).
\]
Using H\"older's inequality we have 
\[
\left|R_{t_{n}}(f)\right|\leq\EE\left(\left[\int_{0}^{t_{n}}(\overline{\tau}(s)-s)~d\Va_{s}(\Xa_{s})~\right]^{p}\right)^{1/p}~\EE\left(\left[e^{\int_{0}^{t_{n}}\Va_{s}(\Xa_{s})ds}+e^{\int_{0}^{t_{n}}\Va_{\underline{\tau}(s)}(\Xa_{\underline{\tau}(s)})ds}\right]^{p^{\prime}}\right)^{1/p^{\prime}}.
\]
To estimate the first term in the r.h.s. of the above estimate, we
use the inequality 
\[
\begin{array}{l}
\EE\left(\left[\int_{0}^{t_{n}}(\overline{\tau}(s)-s)~d\Va_{s}(\Xa_{s})~\right]^{p}\right)^{1/p}\\
\\
\leq\EE\left(\left[\int_{0}^{t_{n}}(\overline{\tau}(s)-s)~\left(\frac{\partial}{\partial s}+L_{s}\right)(\Va_{s})(\Xa_{s})ds~\right]^{p}\right)^{1/p}+\EE\left(\left[\int_{0}^{t_{n}}(\overline{\tau}(s)-s)~dM_{s}(\Va)~\right]^{p}\right)^{1/p}.
\end{array}
\]
Using the generalized Minkowski inequality we prove that 
\[
\begin{array}{l}
\EE\left(\left[\int_{0}^{t_{n}}(\overline{\tau}(s)-s)~\left(\frac{\partial}{\partial s}+L_{s}\right)(\Va_{s})(\Xa_{s})~ds~\right]^{p}\right)^{1/p}\\
\\
\leq\int_{0}^{t_{n}}(\overline{\tau}(s)-s)~\EE\left(\left[\left(\frac{\partial}{\partial s}+L_{s}\right)(\Va_{s})(\Xa_{s})\right]^{p}\right)^{1/p}~ds~\\
\\
\leq\frac{1}{m}~\int_{0}^{t_{n}}\EE\left(\left[\left(\frac{\partial}{\partial s}+L_{s}\right)(\Va_{s})(\Xa_{s})\right]^{p}\right)^{1/p}~ds.
\end{array}
\]
On the other hand, by the Burkholder-Davis-Gundy inequality, for any
$p>0$ we have that 
\begin{eqnarray*}
\EE\left(\left[\int_{0}^{t_{n}}(\overline{\tau}(s)-s)~dM_{s}(\Va)~\right]^{p}\right)^{1/p} & \leq & c_{p}~\EE\left(\left[\int_{0}^{t_{n}}(\overline{\tau}(s)-s)^{2}~d\left[M(\Va)\right]_{s}~\right]^{p/2}\right)^{1/p}\\
 & \leq & \frac{c_{p}}{m}~\EE\left(\left[M(\Va)\right]_{t_{n}}^{p/2}\right)^{1/p}
\end{eqnarray*}
for some finite constant $c_{p}<\infty$ whose values only depends
on the parameter $p$. The end of the proof of the first assertion.
The second one is proved using the decomposition 
\[
\left[\QQ_{t_{n}}^{(m)}-\QQ_{t_{n}}\right](f)=\frac{1}{\Lambda_{t_{n}}^{(m)}(1)}~\left[\Lambda_{t_{n}}^{(m)}\left(f_{n}\right)-\Lambda_{t_{n}}\left(f_{n}\right)\right]
\]
with the centered function $f_{n}=(f-\Lambda_{t_{n}}(f))$. 
\[
\left|\QQ_{t_{n}}^{(m)}(f)-\QQ_{t_{n}}(f)\right|\leq\frac{2}{\Lambda_{t_{n}}^{(m)}(1)\vee\Lambda_{t_{n}}(1)}~\sup_{f\in\Ba_{b}(E_{n})~:~\|f\|\leq1}{\left|\Lambda_{t_{n}}^{(m)}\left(f\right)-\Lambda_{t_{n}}\left(f\right)\right|}.
\]
From these estimates, we find that 
\[
{\left\Vert \overline{r}_{m,t_{n}}\right\Vert _{\rm tv}}\leq\frac{2}{\Lambda_{t_{n}}^{(m)}(1)\vee\Lambda_{t_{n}}(1)}~\left\Vert r_{m,t_{n}}\right\Vert _{\rm tv}
\]
and 
\[
\sup_{m\geq1}{\left\Vert r_{m,t_{n}}\right\Vert _{\rm tv}}\leq a_{p}~b_{t_{n}}^{(p^{\prime})}\left(c_{t_{n}}^{(p)}+d_{t_{n}}^{(p)}\right)
\]
with for any $\frac{1}{p}+\frac{1}{p^{\prime}}=1$ 
\[
b_{t}^{(p^{\prime})}=\EE\left(e^{p^{\prime}\int_{0}^{t}\Va_{s}(\Xa_{s})ds}\right)^{1/p^{\prime}}\vee\EE\left(e^{p^{\prime}\int_{0}^{t}\Va_{\underline{\tau}(s)}(\Xa_{\underline{\tau}(s)})ds}\right)^{1/p^{\prime}}<\infty
\]
as well as 
\[
c_{t}^{(p)}:=\int_{0}^{t}\EE\left(\left[\left(\frac{\partial}{\partial s}+L_{s}\right)(\Va_{s})(\Xa_{s})\right]^{p}\right)^{1/p}ds<\infty\quad\mbox{and}\quad d_{t}^{(p)}:=\EE\left(\left[M(\Va)\right]_{t}^{p/2}\right)^{1/p}<\infty.
\]

Whenever $\Va\in\Ca^{1}([0,\infty[,\Da(L))$, we have 
\[
b_{t}^{(p^{\prime})}\leq\exp{\left(t\sup_{0\leq s\leq t}{\|\Va_{s}\|}\right)}\quad\mbox{{\rm and}}\quad c_{t}^{(p)}\leq t~\sup_{0\leq s\leq t}{\left\Vert \left(\frac{\partial}{\partial s}+L_{s}\right)(\Va_{s})\right\Vert }.
\]
Notice that for stochastic processes with continuous paths, the constant
$d_{t}^{(p)}$ is given by 
\[
d_{t}^{(p)}:=\EE\left(\langle M(\Va)\rangle_{t}^{p/2}\right)^{1/p}=\EE\left(\left(\int_{0}^{t}\Gamma_{L_{s}}(\Va_{s},\Va_{s})(\Xa_{s})ds\right)^{p/2}\right)^{1/p}\leq\sqrt{t}~\sup_{0\leq s\leq t}{\left\Vert \Gamma_{L_{s}}(\Va_{s},\Va_{s})\right\Vert ^{1/2}}.
\]
This ends the proof of the theorem.\cqfd

\subsection{Proof of proposition~\ref{lem-ref-D}}

\label{proof-lem-DC}

In case \textbf{(D)}, the mutation transition only occurs at integer
times. More formally, we have $\Xa_{t_{n}}=X_{\lfloor t_{n}\rfloor}$,
so that for any 
\[
n=km+r\quad\mbox{{\rm with}}\quad r+1<m
\]
we have 
\[
t_{n}=k+r/m\Rightarrow\Ma_{t_{n},t_{n+1}}(y,dz)=\Ma_{k+r/m,k+(r+1)/m}(y,dz)=\delta_{y}(dz).
\]
On the other hand, we have that 
\[
n=km+m-1\Rightarrow t_{n}=(k+1)-1/m\quad\mbox{{\rm with}}\quad t_{n+1}=(k+1)
\]
and 
\[
\Ma_{t_{n},t_{n+1}}(y,dz)=\Ma_{(k+1)-1/m,(k+1)}(y,dz)=M_{k+1}(y,dz).
\]

This ends the proof of the first assertion. Now, we come to the proof
of (\ref{1st-dec-diff}). Firstly, we recall that 
\[
f(t_{n+1},X_{t_{n+1}})=f(t_{n},X_{t_{n}})+\int_{t_{n}}^{t_{n+1}}~\left(\frac{\partial}{\partial s}+L_{s}\right)(f)(s,X_{s})~ds+M_{t_{n+1}}(f)-M_{t_{n}}(f)
\]
for any $f\in\Ca^{1}([t_{n},t_{n+1}],D(L))$, with some martingale
$M_{t}(f)$. This implies that 
\[
\EE_{t_{n},x}\left(f(t_{n+1},X_{t_{n+1}})\right)=f(t_{n},x)+\int_{t_{n}}^{t_{n+1}}~\EE_{t_{n},x}\left[\left(\frac{\partial}{\partial s}+L_{s}\right)(f)(s,X_{s})\right]~ds
\]
where $\EE_{t_{n},x}(\point)$ stands for the conditional expectation
operator given that $X_{t_{n}}=x$. Iterating this formula, we find
that 
\[
\EE_{t_{n},x}\left[\left(\frac{\partial}{\partial s}+L_{s}\right)(f)(s,X_{s})\right]=\left(\frac{\partial}{\partial t_{n}}+L_{t_{n}}\right)(f)(t_{n},x)+\int_{t_{n}}^{s}~\EE_{t_{n},x}\left[\left(\frac{\partial}{\partial r}+L_{r}\right)^{2}(f)(r,X_{r})\right]~dr
\]
as soon as $\left(\frac{\partial}{\partial t}+L_{t}\right)(f)\in\Ca^{1}([t_{n},t_{n+1}],D(L))$.
Under this condition, we find the first order decomposition 
\[
\EE_{t_{n},x}\left(f(t_{n+1},X_{t_{n+1}})\right)=f(t_{n},x)+\left(\frac{\partial}{\partial t_{n}}+L_{t_{n}}\right)(f)(t_{n},x)~\frac{1}{m}+R_{t_{n}}(f)~\frac{1}{m^{2}}
\]
with some remainder operator 
\[
R_{t_{n}}(f)=m^{2}\int_{t_{n}}^{t_{n+1}}\int_{t_{n}}^{s}~\EE_{t_{n},x}\left[\left(\frac{\partial}{\partial r}+L_{r}\right)^{2}(f)(r,X_{r})\right]~dr
\]
such that 
\[
\|R_{t_{n}}(f)\|\leq\sup_{x\in E}\sup_{t_{n}\leq t<t_{n+1}}{\left|\EE_{t_{n},x}\left[\left(\frac{\partial}{\partial t}+L_{t}\right)^{2}(f)(t,X_{t})\right]\right|}.
\]
The end of the proof of (\ref{1st-dec-diff}) is now clear. Using
lemma~\ref{prop-dec}, and the first order decompositions (\ref{1st-exp})
and (\ref{1st-dec-diff}) we have 
\[
L_{t_{n},\mu}^{(m)}=L_{t_{n},\mu}~ \frac{1}{m}+R_{t_{n},\mu}~\frac{1}{m^{2}}
\]
with 
\[
R_{t_{n},\mu}=R_{t_{n}}+\widehat{R}_{t_{n},\mu}+m^{2}\left[\widehat{L}_{t_{n},\mu}~\frac{1}{m}+\widehat{R}_{t_{n},\mu}~\frac{1}{m^{2}}\right]\left[L_{t_{n}}~ \frac{1}{m}+R_{t_{n}}~\frac{1}{m^{2}}\right].
\]
On the other hand, we have the estimates 
\[
m^{2}\left[\widehat{L}_{t_{n},\mu}~\frac{1}{m}+\widehat{R}_{t_{n},\mu}~\frac{1}{m^{2}}\right]\left[L_{t_{n}}~ \frac{1}{m}+R_{t_{n}}~\frac{1}{m^{2}}\right]=\widehat{L}_{t_{n},\mu}\left[L_{t_{n}}+R_{t_{n}}~\frac{1}{m}\right]+\widehat{R}_{t_{n},\mu}L_{t_{n}}^{(m)}
\]
with 
\[
\left\Vert \widehat{R}_{t_{n},\mu}L_{t_{n}}^{(m)}(f)\right\Vert \leq c_{t_{n}}~\left\Vert f\right\Vert 
\]
and 
\begin{eqnarray*}
\left\Vert \widehat{L}_{t_{n},\mu}\left[L_{t_{n}}+R_{t_{n}}~\frac{1}{m}\right](f)\right\Vert  & \leq & c_{t_{n}}\left\Vert L_{t_{n}}(f)+R_{t_{n}}(f)~\frac{1}{m}\right\Vert \\
 & \leq & c_{t_{n}}~\sup_{t_{n}\leq t\leq t_{n+1}}\left(\left\Vert {\partial L_{t}(f)}/{\partial t}\right\Vert +\left\Vert L_{t}(f)\right\Vert +\left\Vert L_{t}^{2}(f)\right\Vert \right)
\end{eqnarray*}
for some finite constant $c_{t_{n}}<\infty$. This ends the proof
of (\ref{1st-dec-diff-C}).

The proof of the last assertion is based on the decomposition 
\[
\Ka_{t_{n},t_{n+1},\mu}\left(\left[f-\Ka_{t_{n},t_{n+1},\mu}(f)(x)\right]^{2}\right)(x)=\Gamma_{L_{t_{n},\mu}^{(m)}}(f,f)(x)-\left(L_{t_{n},\mu}^{(m)}(f)(x)\right)^{2}
\]
and the fact that 
\[
\Gamma_{L_{t_{n},\mu}^{(m)}}(f,f)(x)=\Gamma_{L_{t_{n},\mu}}(f,f)(x)~ \frac{1}{m}+R_{t_{n},\mu}\left((f-f(x))^{2}\right)(x)~~\frac{1}{m^{2}}.
\]
This ends the proof of the proposition.

\cqfd

\subsection{Proof of the first order decompositions (\ref{1st-exp-ref})}

\label{decomp-a}
\begin{itemize}
\item \textbf{Case 1 :} We assume that $\Va_{t}=-\Ua_{t}$, for some non
negative and bounded function $\Ua_{t}$. In this situation, (\ref{MT})
is satisfied by the Markov transition 
\[
\Sa_{t_{n},\mu_{t_{n}}^{(m)}}(x,dy):=e^{-\Ua_{t_{n}}(x)/m}~\delta_{x}(dy)+\left(1-e^{-\Ua_{t_{n}}(x)/m}\right)~\Psi_{e^{-\Ua_{t_{n}}/m}}(\mu_{t_{n}}^{(m)})(dy).
\]
In this situation, we notice that 
\begin{eqnarray*}
m~\left[\Sa_{t_{n},\mu_{t_{n}}^{(m)}}(f)(x)-f(x)\right] & = & m~\left(1-e^{-\Ua_{t_{n}}(x)/m}\right)~\left[\Psi_{e^{-\Ua_{t_{n}}/m}}(\mu_{t_{n}}^{(m)})(f)-f(x)\right]\\
 & = & \Ua_{t_{n}}(x)~\int~\left(f(y)-f(x)\right)~\mu_{t_{n}}^{(m)}(dy)+O(1)=\widehat{L}_{t_{n},\mu_{t_{n}}^{(m)}}(f)+O(1).
\end{eqnarray*}
Much more is true; if we set 
\begin{eqnarray*}
\widehat{R}_{t_{n},\mu_{t_{n}}^{(m)}}(f) & := & m\left[m~\left[\Sa_{t_{n},\mu_{t_{n}}^{(m)}}(f)-f\right]-\widehat{L}_{t_{n},\mu_{t_{n}}^{(m)}}(f)\right]\\
 & = & m\left[m~\left(1-e^{-\Ua_{t_{n}}/m}\right)~\left[\Psi_{e^{-\Ua_{t_{n}}/m}}(\mu_{t_{n}}^{(m)})(f)-f\right]-\Ua_{t_{n}}\left[\mu_{t_{n}}^{(m)}(f)-f\right]\right]
\end{eqnarray*}
then we find that 
\[
\begin{array}{l}
\left|\widehat{R}_{t_{n},\mu_{t_{n}}^{(m)}}(f)\right|\\
\\
\leq m\left|m~\left(1-e^{-\Ua_{t_{n}}/m}\right)-\Ua_{t_{n}}\right|\times\left|\Psi_{e^{-\Ua_{t_{n}}/m}}(\mu_{t_{n}}^{(m)})(f)-f\right|+\Ua_{t_{n}}~m\left|[\Psi_{e^{-\Ua_{t_{n}}/m}}(\mu_{t_{n}}^{(m)})(f)-\mu_{t_{n}}^{(m)}(f)\right|
\end{array}
\]
Using the fact that 
\[
m~\left[\Psi_{e^{-\Ua_{t_{n}}/m}}(\mu_{t_{n}}^{(m)})-\mu_{t_{n}}^{(m)}\right](f)=\frac{1}{\mu_{t_{n}}^{(m)}\left(e^{-\Ua_{t_{n}}/m}\right)}~\mu_{t_{n}}^{(m)}\left(m~\left[e^{-\Ua_{t_{n}}/m}-1\right]\left[f-\mu_{t_{n}}^{(m)}(f)\right]\right)
\]
we prove the following first order expansion 
\begin{equation}
\Sa_{t_{n},\mu_{t_{n}}^{(m)}}(f)-f=\widehat{L}_{t_{n},\mu_{t_{n}}^{(m)}}(f)~\frac{1}{m}+\widehat{R}_{t_{n},\mu_{t_{n}}^{(m)}}(f)~\frac{1}{m^{2}}\label{1st-exp}
\end{equation}
with some integral operator $\widehat{R}_{t_{n},\mu_{t_{n}}^{(m)}}$
such that 
\[
\sup_{m\geq1}{\left\Vert \widehat{R}_{t_{n},\mu_{t_{n}}^{(m)}}(f)\right\Vert }\leq c~\left\Vert \Ua_{t_{n}}\right\Vert ^{2}\mbox{{\rm osc}}(f).
\]

\item \textbf{Case 2 :} We assume that $\Va_{t}$ is a positive and bounded
function. In this situation, (\ref{MT}) is satisfied by the Markov
transition 
\[
\Sa_{t_{n},\mu_{t_{n}}^{(m)}}(x,dy):=\frac{1}{\mu_{t_{n}}^{(m)}\left(e^{\Va_{t_{n}}/m}\right)}~\delta_{x}(dy)+\left(1-\frac{1}{\mu_{t_{n}}^{(m)}\left(e^{\Va_{t_{n}}/m}\right)}\right)~\Psi_{\left(e^{\Va_{t_{n}}/m}-1\right)}(\mu_{t_{n}}^{(m)})(dy).
\]
In this situation, we notice that 
\begin{eqnarray*}
m~\left[\Sa_{t_{n},\mu_{t_{n}}^{(m)}}(f)(x)-f(x)\right] & = & m~\left(1-\frac{1}{\mu_{t_{n}}^{(m)}\left(e^{\Va_{t_{n}}/m}\right)}\right)~\left[\Psi_{\left(e^{\Va_{t_{n}}/m}-1\right)}(\mu_{t_{n}}^{(m)})(f)-f(x)\right]\\
 & = & \int~\left(f(y)-f(x)\right)~\Va_{t_{n}}(y)~\mu_{t_{n}}^{(m)}(dy)+O(1).
\end{eqnarray*}
Using some elementary calculations, we also prove a first order expansion
of the same form as in (\ref{1st-exp}). 
\item \textbf{Case 3:} The Markov transport equation (\ref{MT}) is satisfied
by the transitions 
\[
\begin{array}{l}
\Sa_{t_{n},\mu_{t_{n}}^{(m)}}(x,dy):=\left(1-\frac{\mu_{t_{n}}^{(m)}\left(\left(e^{\Va_{t_{n}}/m}-e^{\Va_{t_{n}}(x)/m}\right)_{+}\right)}{\mu_{t_{n}}^{(m)}\left(e^{\Va_{t_{n}}/m}\right)}\right)~\delta_{x}(dy)\\
\\
\hskip4cm+{\displaystyle \frac{\mu_{t_{n}}^{(m)}\left(\left(e^{\Va_{t_{n}}/m}-e^{\Va_{t_{n}}(x)/m}\right)_{+}\right)}{\mu_{t_{n}}^{(m)}\left(e^{\Va_{t_{n}}/m}\right)}~\Psi_{\left(e^{\Va_{t_{n}}/m}-e^{\Va_{t_{n}}(x)/m}\right)_{+}}(\mu_{t_{n}}^{(m)})(dy)}.
\end{array}
\]

The above Markov transition is well defined since we have 
\[
\frac{\mu_{t_{n}}^{(m)}\left(\left(e^{\Va_{t_{n}}/m}-e^{\Va_{t_{n}}(x)/m}\right)~1_{\Va_{t_{n}}>\Va_{t_{n}}(x)}\right)}{\mu_{t_{n}}^{(m)}\left(e^{\Va_{t_{n}}/m}1_{\Va_{t_{n}}>\Va_{t_{n}}(x)}\right)+\mu_{t_{n}}^{(m)}\left(e^{\Va_{t_{n}}/m}1_{\Va_{t_{n}}\leq\Va_{t_{n}}(x)}\right)}\leq1
\]
as soon as 
\[
\mu_{t_{n}}^{(m)}\left(e^{\Va_{t_{n}}/m}1_{\Va_{t_{n}}\leq\Va_{t_{n}}(x)}\right)+\mu_{t_{n}}^{(m)}\left(e^{\Va_{t_{n}}(x)/m}~1_{\Va_{t_{n}}>\Va_{t_{n}}(x)}\right)=\mu_{t_{n}}^{(m)}\left(e^{[\Va_{t_{n}}\wedge\Va_{t_{n}}(x)]/m}\right)>0.
\]
Also notice that 
\[
\mu_{t_{n}}^{(m)}\left(\left(e^{\Va_{t_{n}}/m}-e^{\Va_{t_{n}}(x)/m}\right)_{+}\right)=\mu_{t_{n}}^{(m)}\left(\Va_{t_{n}}>\Va_{t_{n}}(x)\right)\times\left[\mu_{t_{n}}^{(m)}\left(e^{\Va_{t_{n}}/m}~|~\Va_{t_{n}}>\Va_{t_{n}}(x)\right)-e^{\Va_{t_{n}}(x)/m}\right]
\]
so that a particle in some state $x$ is more likely to be recycled
when the potential values $\Va_{t_{n}}(\overline{\Xa}_{t_{n}})$ of
random states $\overline{\Xa}_{t_{n}}$ with law $\mu_{t_{n}}^{(m)}$
are more likely to be larger that $\Va_{t_{n}}(x)$.

Finally, for bounded potential functions we observe that 
\[
\begin{array}{l}
m~\left[\Sa_{t_{n},\mu_{t_{n}}^{(m)}}(f)(x)-f(x)\right]\\
\\
=m~{\displaystyle \frac{\mu_{t_{n}}^{(m)}\left(\left(e^{\Va_{t_{n}}/m}-e^{\Va_{t_{n}}(x)/m}\right)_{+}\right)}{\mu_{t_{n}}^{(m)}\left(e^{\Va_{t_{n}}/m}\right)}~\left[\Psi_{\left(e^{\Va_{t_{n}}/m}-e^{\Va_{t_{n}}(x)/m}\right)_{+}}(\mu_{t_{n}}^{(m)})(f)-f(x)\right]}\\
\\
={\displaystyle \int~[f(y)-f(x)]~(\Va_{t_{n}}(y)-\Va_{t_{n}}(x))_{+}~\mu_{t_{n}}^{(m)}(dy)+O(1)}.
\end{array}
\]
Using some elementary calculations, we also prove a first order expansion
of the same form as in (\ref{1st-exp}). 

\end{itemize}
\bibliographystyle{plain}
\bibliography{discrete-continuous}

\end{document}